\documentclass[12pt,draftcls,onecolumn]{IEEEtran}
\usepackage{amsmath,amsbsy,amscd,amssymb}
\usepackage[pdftex]{graphicx}
\newtheorem{theorem}{Theorem}

\newtheorem{example}{Example}
\newtheorem{lemma}{Lemma}

\newcommand{\dessin}[4]{
\begin{figure}[htb]
\begin{center}
\includegraphics[scale=#2]{#1}
\caption{#3}
\label{#4}
\end{center}
\end{figure}}
\newcommand{\mpdiv}{/}

\newcommand{\CC}{{\mathcal C}}

\newcommand{\UU}{{\mathcal U}}
\newcommand{\VV}{{\mathcal V}}

\newcommand{\YY}{{\mathcal Y}}
\newcommand{\QQ}{{\mathcal Q}}

\newcommand{\PP}{{\mathcal P}}

\newcommand{\GG}{{\mathcal G}}

\newcommand{\ZZ}{{\mathcal Z}}

\newcommand{\bbR}{{\mathbb R}}

\newcommand{\bbN}{{\mathbb N}}

\newcommand{\bbT}{{\mathbb T}}

\newcommand{\Rbmin}{\overline{\bbR}_{\min}}

\newcommand{\Rmin}{\bbR_{\min}}

\def\c01{\CC([0,1])}

\graphicspath{{figures/}}

\begin{document}

\title{About Dynamical Systems Appearing in the Microscopic Traffic Modeling}
\author{Nadir Farhi, Maurice Goursat \& Jean-Pierre Quadrat\\
\small{INRIA - Paris - Rocqencourt}\\
Domaine de Voluceau, 78153, Le Chesnay, Cedex France.\\
\texttt{nadir.farhi@inria.fr}}
\date{arXiv: January 26$^{\text{th}}$, 2010}
\maketitle

\begin{abstract}
Motivated by microscopic traffic modeling, 
we analyze dynamical systems which have a
piecewise linear concave dynamics not necessarily monotonic. We
introduce a deterministic Petri net extension where edges may have
negative weights. The dynamics of these Petri nets are well-defined
and may be described by a generalized matrix with a submatrix in
the standard algebra with possibly negative entries, and another
submatrix in the minplus algebra. When the dynamics is additively
homogeneous, a generalized additive eigenvalue may be introduced,
and the ergodic theory may be used to define a growth rate under
additional technical assumptions. In the traffic example of two
roads with one junction, we compute explicitly the eigenvalue 
and we show, by numerical simulations,
that these two quantities (the additive eigenvalue and the growth rate) 
are not equal, but are close to each other.
With this result, we are able to extend the well-studied notion of
fundamental traffic diagram (the average flow as a function of the
car density on a road) to the case of two roads with one junction and
give a very simple analytic approximation of this diagram where
four phases appear with clear traffic interpretations.
Simulations show that the fundamental diagram shape obtained is also valid 
for systems with many junctions. To simulate these systems,
we have to compute their dynamics, which are not quite simple. 
For building them in a modular way,
we introduce generalized parallel, series and feedback
compositions of piecewise linear concave dynamics.
\end{abstract}

\section{Introduction}

The main purpose of this paper is to explore some new points on
dynamical systems that have emerged with the study of microscopic
traffic modeling using minplus algebra and Petri nets. These
points are the following~:
\begin{itemize}
\item Standard Petri nets or linear maxplus algebra have not
enough modeling power to describe the dynamics of a  traffic
system as simple as two roads with a crossing and with the simplest
vehicle move. \item A remedy to this difficulty is to introduce
possibly negative weights on Petri net arcs (not to be confused
with negative places or negative tokens as in~\cite{LeePark}) or
to use nonlinear minplus dynamics that are not monotonic. \item A
consequence is that all known standard methodology based mainly on
dynamic programming is ineffective. \item Nevertheless, it is
sometimes possible to use the standard ergodic theory, and/or to
compute a generalized eigenvalue for the subclass of additively
homogeneous systems and get nice qualitative results on the
system. \item When it is not possible to derive analytical
results, we can do numerical simulations. In general, a
microscopic modeling of a realistic system has a very large number
of state variables difficult to manage. However, for the
construction of the model, a generalized matrix calculus using
minplus algebra can be developed for the class of system described
here (having piecewise affine concave dynamics).
\end{itemize}
All these points are illustrated on the traffic application for
which we have obtained a new result: a good analytic approximation
of the macroscopic law called the fundamental diagram in the case of
roads with junctions.

The traffic on  a road has been studied from different points of
view on a macroscopic level. For example:
\begin{itemize}
\item  The  Lighthill-Whitham-Richards  Model~\cite{LW}, which is the most
standard view, expresses the mass conservation of cars seen as a
fluid:
$$
\begin{cases}
\partial_t \rho+\partial_x \varphi=0\;, \\ \varphi=f(\rho),
\end{cases}
$$ where $\varphi(x,t)$ denotes the  flow at time $t$ and position $x$
on  the road;  $\rho(x,t)$ denotes  the density;  and $f$  is  a
given function  called  the \emph{fundamental  traffic  diagram}.  For
traffic, this diagram plays the role of the perfect gas law for
the fluid dynamics. The diagram has been estimated using
experimental data, and its behavior is quite different from
standard gas at high density. We obtain an idea of the diagram's
shape in the subsection titled ``traffic on a circular road.''

Daganzo in \cite{Dag05} using the variable $N(t,x)$ counting the cumulated
number of vehicles having reached the point $x$ at time $t$, and remarking
that $\varphi=\partial_t N$ and $\rho=-\partial_x N$ has interpreted
the equation $\varphi=f(\rho)$ as an Hamilton-Jacobi equation  when $f$
is concave. The other equation $\partial_t \rho+\partial_x \varphi=0$
telling that $\partial_{tx} N=\partial_{xt}N$. 
Then, in \cite{Dag08}, Daganzo and Geroliminis generalized this to a network by postulating
the existence of an Hamiltonian in this case also and by trying to approximate
it as an infimum of linear constraints that he obtains by physical considerations
on the street involved. In \cite{GD07,GD08}, Geroliminis and Daganzo studied what happen
experimentally on real towns. 

\item  The  kinetic   model  (Prigogine-Herman~\cite{PRI})  gives
the evolution of the  density of particles $\rho(t,x,v)$ as  a
function of $t,x$ and $v$, where $v$ is the speed of particles.
The model is given by:
$$\partial_t    \rho+v\partial_x    \rho=C(\rho,\rho)\;   ,$$    where
$C(\rho,\rho)$ is an interacting term that is, in general,
quadratic in $\rho$ and, as such, models the driver behavior in a
simple way. From the distribution $\rho$, we can derive all the
useful quantities, such as the average speed $\bar{v}(t,x)=\int
v\rho(t,x,v)dv$.
\end{itemize}
The  integration of the second model  is  more time consuming and,
therefore, not  used in  practice. The main
interest of the Prigogine-Herman equation is that we can derive
macroscopic laws like the fundamental diagram from its solution.

The first  model (the Lighthill-Whitham-Richards  Model~\cite{LW}) assumes the
knowledge of the  fundamental traffic diagram  $f$.  
This function has been studied
not only experimentally but also theoretically using simple  microscopic models
such as exclusion processes~\cite{DER,  BLA},  cellular automata
\cite{CL,CSS}, or simulation of individual  car dynamics. 
Here, we first recall a way to derive an approximation of this diagram on a circular road,
then we present an extension to the two dimensional traffic (crossing roads).
This derivation consists of computing the eigenvalue of a simple
minplus linear  system counting the number  of vehicles entering
in a road section.

The main result obtained is the generalization of this fundamental
law to  the 2D  case  where the  roads  cross each  other. In
statistical physics, a lot of numerical work
(\cite{Cues,Cues2,Cues3,BARL,BBSS}) and good
surveys~\cite{CSS,HEL} have been done on idealized towns. These
works analyze numerical experimentations based on various
stochastic models with or without turning possibilities and show
the existence of a threshold of the density at which the system
blocks suddenly. The particular case of one junction is studied
in~\cite{FUKUI,FUKUI1} where precise results are given for the
stochastic case without turning possibilities. Let us mention
also the attempt to derive a 2D fundamental diagram by Helbing~\cite{Hel09}
using queuing theory. In the saturated and unsaturated traffic case he derives
a formula for an intersection that he extends to an area of a town
that fits well to experimental datas. 

Here we present a deterministic model with turning possibilities, based on Petri
nets and minplus algebra. The minplus linear model on a unique
road can be described in terms of event graph which is a subclass
of Petri nets. The presence of conflicts at junctions prevents the extension of
this model to the 2D case. 

We propose a way to solve the difficulty by extending the class of weights 
used in Petri nets allowing negative weights. Due to
such weights, the firing of a transition can consume tokens
downstream, and the modeling power of the Petri nets is improved
significantly. This possibility is used to model the authorization
to enter in the junctions. 

The dynamics  of general Petri nets
allowing negative weights can be written easily but is neither
linear in minplus algebra nor monotone. Nevertheless, in the
traffic applications given here, dynamics ($x^{k+1}=f(x^k)$) are
always homogeneous of degree one ($f(\lambda \otimes x)=\lambda
\otimes f(x)$, where $\otimes$ denotes the minplus multiplication
that is the standard addition). For such systems, the eigenvalue
problem (computation of $\lambda$ such that $\lambda \otimes
x=f(x)$) can be reduced to a fixed point problem. The existence of
the growth rate $\chi=\lim_k x^k/k$ is due to the existence of a
Birkhoff average. The quantities $\lambda$ and $\chi$ coincide
when $f$ is also monotone, but this is not true in the general
case. The monotone case has been studied carefully
in~\cite{GG04,MPN}. In traffic examples of roads with junctions,
the dynamics is not monotone. 


In all the traffic examples, we study systems of roads 
on a torus without entries, in such a way
that the number of vehicles remains constant in the system. Thus,
we represent an idealization of constant densities for a
more realistic system. To maintain a constant density in an open
system, a new vehicle has to enter each time another one leaves
the system. This is mathematically equivalent to consider
circular roads. If we want to maintain a constant local density in
a large system, it is the same problem; each time a vehicle leaves
the local subsystem, another one has to enter.

We study indetail the particular case of two circular roads with one junction.
For this system, a result on the existence of the
growth rate is obtained using the nondecreasing trajectory property
of the states starting from zero. From this result, we show that
the distances between the states stay bounded. The eigenvalue problem
can be solved explicitly. On simulation, we see that the eigenvalue
and the growth rate generally do not coincide.  However, these
quantities are very close for any fixed density. 
Therefore, the simple formulas obtained
for the eigenvalue give a good approximation of the traffic
fundamental diagram. 

This definition of the fundamental diagram must not be confused with the 
one used by Daganzo~\cite{Dag05}, which is the Hamiltonian of the traffic dynamics
interpreted as a dynamic programming equation in the simplest case of a single road
and generalized in heuristic way to the city case.
His fundamental diagram mainly coincides with our minplus dynamics,
in the case of one street, but is observed experimentally 
not derived of a microscopic model as we do here.
The dependence of the eigenvalue 
or the average growth rate with the density of vehicles in the system,
which is a difficult result to obtain even in the case of two streets
with one intersection, is not done in the Daganzo work.
Moreover, for the city case, we remark also, on the simple example studied here, 
that the dynamics is not a standard dynamic programming equation associated to a deterministic 
or a stochastic control problem (since the dynamics is not monotone) which could be 
a problem, at least at the conceptual level, for the point of view adopted by Daganzo.

The fundamental diagram that we obtain for a system of two roads with one junction, shows four phases:
\begin{itemize}
\item \emph{Free phase}: The density is low where the vehicles do not
interact.
\item \emph{Saturation phase}: The junction is saturated, but the locations downstream
of the junction are freed when a vehicle wants to leave the junction.
\item \emph{Recession phase}: The locations downstream of the junction are
sometimes crowded when a vehicle wants to leave the junction.
\item \emph{Freeze phase}: The vehicles cannot move.
\end{itemize}
These phases are derived from a unique simple model and are not postulated to
obtain different models subsequently analyzed.

Preliminary results on the traffic dynamics have been presented
in~\cite{FAR,FAR2,FAR3}. In \cite{FARTh}, many developments,
complementary results, and other completely solved examples can be
found. The theorems given here on the eigenvalue problem and the
growth rate complete some of the main results of~\cite{FARTh} by
relaxing some assumptions and clarifying the growth rate existence
problem.

We show that a system theory can be developed  for the class of
concave polyhedral 1-homogeneous systems, which we construct using
parallel, series, and feedback compositions. Moreover, the
dynamics of these systems are characterized by the composition of
a standard matrix and a minplus matrix. The sum of the entries in
each line of the standard matrix is equal to one. We can compute
easily the transformation on these two matrices corresponding to
the three composition operations. The interest of this ``system
theory'' can be shown by building the traffic dynamics of a
regular town starting from three elementary systems. We do not
detail the construction, but we show a stationary vehicle
distribution obtained for a simple regular town with dynamics
built in that way. In this case, the fundamental diagrams still
present the four phases observed in the simple case of two roads
with one junction. All the numerical simulations have been done
using the ScicosLab software~\cite{ScicosLab}.

This paper is divided into three parts. First, we recall the basic
results of minplus algebra; we present a system theory for
polyhedral concave 1-homogeneous systems and discuss their growth
rate and eigenvalue problems. Second, we present Petri nets
allowing negative weights. Third, we give applications to the
computation of the fundamental traffic diagram; we consider the
case of one circular road with minplus linear dynamics and study
the case of two circular roads with a junction in detail; then, we
explain briefly the way to build the dynamics using the system
theory described in the first part.

\section{Minplus Algebra and extensions}
\subsection{Review of minplus algebra}
In this section, we revisit the main definitions and properties  of
the minplus  algebra.  An  in-depth  treatment of  the subject is in
\cite{BCOQ}.

The structure \( \bbR_{\min }=\left( \bbR\cup \left\{ +\infty
\right\} ,\oplus ,\otimes \right)\) is defined  by the set \( \bbR
\cup \left\{ +\infty \right\} \) endowed with the operations \(
\min \) (denoted by \(  \oplus \),  called minplus  \emph{sum})
and  \( + \)  (denoted by \(\otimes  \),   called  minplus
\emph{product}). The  element  \( \varepsilon =+\infty  \) is the
\emph{zero} element \( \varepsilon  \oplus  x=x \), and is
\emph{absorbent} \(  \varepsilon \otimes x=\varepsilon  \). The
element  \( e=0 \) is  the \emph{unit} element \( e\otimes x=x \).
The main  difference with respect to the conventional algebra is
the idempotency of the addition \( x\oplus x=x  \)  and the fact
that  the  addition  cannot  be \emph{simplified} $a\oplus
b=c\oplus  b  \not \Rightarrow a=c$. This structure is called
\emph{minplus algebra}. We will call $\Rbmin$ the completion of
$\Rmin$ by $-\infty$ with $-\infty\otimes\varepsilon=\varepsilon$.

This  minplus  structure on  scalars  induces  an idempotent
semiring structure on $m\times m$ square matrices with the
element-wise minimum denoted by $\oplus$ and   the    matrix
product defined   by \(   \left(   A\otimes B\right)_{ik}=\min
_{j}\left( A_{ij}+B_{jk}\right)\), where  the zero and the unit
matrices are still denoted by $\varepsilon$ and $e$.  We associate
a precedence graph  $\GG(A)$ to  a square matrix  $A$  where the
nodes of the graph correspond to the columns of  the matrix $A$
and the edges of the graph correspond to the non-zero ($\neq
\varepsilon$) entries of the matrix. The \emph{weight} of the edge
going from $i$ to $j$ is the non-zero entry $A_{ji}$. We define
the weight of a path $p$ in $\GG(A)$, which we denote by
$|p|_{w}$, as the minplus product of the weights of the  edges
composing  the path (that is the standard sum of weights).  The
number of  edges of a path $p$ is denoted  by $|p|_{l}$.   We will
use the following fundamental result (see \cite{BCOQ}).
\begin{theorem}\label{eigenth}
If the graph  $\GG(A)$ associated to an $m\times  m$ minplus matrix
$A$  is  strongly connected,  then  the  matrix  $A$ admits  a  unique
\emph{eigenvalue}  $\lambda\in\bbR_{\min  }\setminus \{\varepsilon\}$:
\begin{equation}
  \exists\   X\in\bbR_{\min   }^{m},   X\neq  \varepsilon:   A\otimes
  X=\lambda      \otimes       X      \;      \text{       with      }
  \lambda=\min_{c\in\CC}\frac{|c|_{w}}{|c|_{l}}\;,
  \label{eigen}
\end{equation}
\normalsize where $\CC$ is the set of circuits of $\GG(A)\;\square $
\end{theorem}

\subsection{A generalized matrix calculus}\label{MaxPlus}
In a Petri net, two kinds of operations appear. One is the
accumulation of resources in the places, and the other is the
synchronization in the transitions. The first operation can be
modeled by addition, and the second by a min or max (a task can
start at the maximum of the arrival instants of the resources
needed by the task). Matrix notations can be generalized to such
situations.

We consider the set of $m\times m$ matrices where the rows [resp. columns] 
are partitioned in two sets the standard and the minplus
(here the $m'$ first  rows [resp. columns], and the last $m"$
rows [resp. columns]) with entries in $\bbR\cup \left\{ +\infty
,-\infty \right\}$, equipped with the two operations $\boxplus$
and $\boxtimes$ defined by:
\[ \begin{bmatrix}A & B \\ C & D\end{bmatrix}\boxplus \begin{bmatrix}A' & B' \\ C' & D'\end{bmatrix}=
\begin{bmatrix}A+A' & B+B' \\ C\oplus C' & D\oplus D'\end{bmatrix} \; ,\]

\[ \begin{bmatrix}A & B \\ C & D\end{bmatrix}\boxtimes \begin{bmatrix}A' & B' \\ C' & D'\end{bmatrix}=
\begin{bmatrix}AA'+BC' & AB'+BD' \\ C\otimes A'\oplus D \otimes C' & C \otimes B'\oplus D \otimes D'\end{bmatrix} \; .\]

Since the entries are in an extension of $\bbR$, we have to specify
the scalar addition and multiplication table:
$$0\times \pm \infty=\pm \infty \times 0=0, \quad +\infty\otimes (-\infty)=+\infty -\infty =+\infty\; .$$
These choices have been made to preserve the absorption properties
for the multiplication of the null elements of the standard
algebra ($0$) and of the minplus algebra ($\varepsilon=+\infty$).
This absorption property is useful to model the absence of an arc
in the precedence graph $\GG(A)$ associated to a square matrix $A$
(defined in the same way as in the pure minplus case).

The addition $\boxplus$ is associative, commutative and has the null element
$\begin{bmatrix} 0 & 0 \\ \varepsilon & \varepsilon \end{bmatrix}$ still denoted $\varepsilon$.

The multiplication $\boxtimes$ has no identity element.
It is neither associative nor commutative nor distributive with respect to the addition.

The main interest of this operation appears when the left operand is a vector 
$Y=A\boxtimes X$, where $X\in \Rbmin^m$ is a vector and $A$
is a $m\times m$ matrix with entries in $\Rbmin$.
the vector $Y$, seen as a function of $X$, is a set of $m'$ standard linear forms 
and of $m"$ minplus linear forms on
$\Rbmin^m$ with $m=m'+m"$.
The operation $Z=A\boxtimes(A\boxtimes X)$  corresponds to the compositions of these linear forms
but these compositions do not define anymore a set of standard and minplus linear forms and
$Z\neq(A\boxtimes A)\boxtimes X$.

Moreover, when the left operand is a vector, this multiplication can be represented 
by the graph $\GG(A)$ where we have
two kinds of nodes (those corresponding to the standard $+$
operation, and those corresponding to $\oplus$ operation) and two
kinds of edges (those which operate multiplicatively ($\times$) and
those which operate additively ($\otimes$)).
\begin{example}\label{2Dexample}
Let us consider the graph $\GG(A)$ associated to the matrix
 $\begin{bmatrix} a &  b \\ c & d \end{bmatrix}$
 (Figure-\ref{Bigraph})
 with one node associated to the standard algebra and one node to the minplus algebra.
 Then $y=A\boxtimes x$, where $y$ and $x$ are two vectors with two entries, means:
$$y_1=ax_1+bx_2,\quad y_2=\min(c+x_1,d+x_2)\; .$$
\dessin{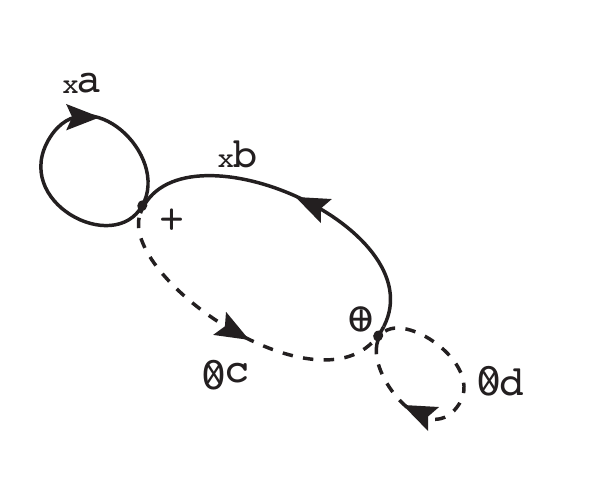}{1}{The incidence graph associated to the matrix $A$ with one $\oplus$node and one +node.}{Bigraph}
\end{example}

We adopt the following conventions to solve some ambiguities in the
formula notations:

-- As soon as a minplus symbol appears in a set of formulas, all
the operations must be understood in the minplus sense with the
exception of the exponent that must be understood in the standard
sense. For example, $x^{a/b}\oplus 1$ or $\sqrt[b]{x^a}\oplus 1$
must be understood as $\min((a/b)x,1)$ and not as $\min((a-b)x,1)$.
In minplus sense, $a/b=a-b$ since it is the solution\footnote{For 
readers familiar residuation theory, the minplus division used
here means the standard minus operator with the convention
previously given for infinite elements. This choice is incompatible
with the residuation which chooses the smallest solution of
$b\otimes x=a$.} of $b\otimes x =a$. The equation $b\otimes x =a$
means $b+x=a$.

-- In a minplus formula, the rational numbers (written with
figures) are denoted in the standard algebra. For example, $\frac{1}{2}x\oplus
1$ (instead of $\sqrt{1} x \oplus 1$) means $\min(0.5+x,1)$ and not
$\min(-1+x,1)$, but $(a/b)x\oplus 1$ means $\min(a-b+x,1)$.

-- A non-zero element in the minplus sense means a finite element in
the usual sense. Positive element will be used always in the
standard sense.

With these conventions the formulas are concise and not ambiguous.

\subsection{A generalized system theory}
We can develop a generalized system theory based on the two
operations $\boxplus$ and $\boxtimes$. We partition the states
[resp. inputs, outputs]  into two classes, the standard states
[resp. inputs, outputs] and the minplus states [resp. inputs,
outputs]. Then the dynamics can be written:
\begin{align*} \begin{bmatrix}
           X^{k+1}\\  Y^{k+1}
   \end{bmatrix}&=
          \begin{bmatrix}
            A &  B\\ C & \varepsilon
          \end{bmatrix}\boxtimes
           \begin{bmatrix}
           X^{k}\\ U^{k}
          \end{bmatrix}
\end{align*}
These dynamics, denoted by $S$ and defined by the matrices
$(A,B,C)$, associate the output signals $(Y^k)_{k\in\mathbb N}$ to
the input signals $(U^k)_{k\in\mathbb N}$. We write $Y=S(U)\;$.

Let us note that $\boxtimes$ corresponds to the definition given
previously when there is a permutation of rows and columns.  This
is because $[X^k\;U^k]'$ [resp. $[X^k\;Y^k]'$)] is not written in
the canonical form since all the standard entries are not followed
by all the minplus entries but by the standard states, the minplus
states, the standard inputs [resp. outputs], and the minplus
inputs [resp. outputs]. See the Petri net
application~(\ref{PetriProd}) for developed formulas.

On these systems, we define the following operations:
\begin{itemize}
  \item  \textbf{Parallel Composition.}  Given  two systems  $S_1$ and
  $S_2$ (with the same numbers of inputs and outputs), we define the system $S=S_{1}\boxplus S_{2}$ obtained by using
  the same entries  and adding\footnote{The  standard output nodes are
  added in  the standard  algebra and the  minplus output nodes are added in the minplus
  algebra.} the outputs. The dynamics of $S$ is defined by:
  $$A=\begin{bmatrix} A_1 & \varepsilon\\ \varepsilon & A_2
      \end{bmatrix},\;
   B=\begin{bmatrix} B_1\\ B_2
           \end{bmatrix},\;
   C=\begin{bmatrix} C_1 & C_2
           \end{bmatrix}\;.$$

  \item  \textbf{Series  Composition.}  Given  two  systems $S_1$  and
  $S_2$ (where the output numbers of $S_2$ equal the input numbers of $S_1$),  we define  the  system $S=S_{1}\boxtimes  S_{2}$ obtained  by
  composition  of   the  two  systems:   $S(U)=S_1(S_2(U))$. Using the new state $[X_1\; X_2\; Y_2]'$, the
  dynamics of $S$ is defined by:
  $$A=\begin{bmatrix} A_1 & \varepsilon & B_1\\ \varepsilon & A_2 & \varepsilon\\ \varepsilon & C_2 & \varepsilon
           \end{bmatrix},\;
   B=\begin{bmatrix} \varepsilon \\ B_2\\ \varepsilon
           \end{bmatrix},\;
   C=\begin{bmatrix} C_1 & \varepsilon & \varepsilon
           \end{bmatrix}\; . $$

  \item \textbf{Feedback.} Given a  system $S$, we define the feedback
  system    $S^{\boxdot}$   as   the    solution   in    $Y$   of
  $Y=S(U\boxplus Y)$.  Using the new state $[X\; Y]'$, the dynamics of $S^{\boxdot}$ are defined by:
  $$A^\boxdot=\begin{bmatrix} A & B\\ C & \varepsilon
           \end{bmatrix},\;
   B^\boxdot=\begin{bmatrix} B\\ \varepsilon
           \end{bmatrix},\;
  C^\boxdot=\begin{bmatrix}C        &        \varepsilon        \end{bmatrix}\;. $$
\end{itemize}
\subsection{Additively homogeneous autonomous dynamic systems}
Let us now discuss one-homogeneous minplus dynamical systems which
are a large class applied frequently when there is a conservation
such as probability mass or tokens in Petri nets. Because
we accept negative entries, we obtain a generalization of two
cases. One case is a measure not necessarily positive with a total
mass equal to one.  The other case is a conservation of tokens in
a situation where negative entries may appear. Let us start with
an academic example. The traffic modeling will give more concrete
examples in the following sections.
\begin{example}
Let us go back to Example \ref{2Dexample} and assuming that
$a+b=1$. Adding $\lambda$ to each component of $x$ (that we can
write $\lambda \otimes x$) implies that the two components of $y$
are augmented by $\lambda$. We have $$A\boxtimes (\lambda \otimes
x)=\lambda \otimes (A\boxtimes x)\;.$$ We say that the system is
\emph{additively homogeneous of degree 1} or, more simply,
\emph{homogeneous}. Indeed, using the minplus notation,
$y=A\boxtimes x$ can be written:
$$y_1=(x_1)^{\otimes a}\otimes (x_2)^{\otimes b},\quad y_2=c\otimes x_1 \oplus d\otimes x_2\; ,$$
which is clearly homogeneous of degree 1 in the minplus algebra as
long as $a+b=1$. Moreover, if $a$ and $b$ are nonnegative, then
the transformation is \emph{nondecreasing}.

To simplify the notations, we will write the transformation in the
following way:
$$y_1=(x_1)^{a}(x_2)^{b},\quad y_2=cx_1 \oplus d x_2\; .$$
\end{example}

More generally, we say that a function $f:\Rbmin^{n} \mapsto
\Rbmin^{n}$ is \emph{homogeneous} if $$f(\lambda \otimes x)=\lambda
\otimes f(x)\;.$$

\subsection{Eigenvalues of homogeneous systems}
The \emph{eigenvalue problem} of a function $f$ can be formulated
as finding non-zero $x\in \Rmin^n$ and $\lambda \in \Rbmin$ such
that:
$$\lambda \otimes x=f(x)\;  .$$ Since $f$ is  homogeneous, we can suppose without
loss of  generality that if  $x$ exists then  $x_{1}\neq
\varepsilon$. The eigenvalue problem becomes:
$$\begin{cases}         \lambda&=f_{1}(x/x_{1})\;         ,         \\
x_{2}/x_{1}&=(f_{2}/f_{1})(x/x_1)\;     ,     \\     \cdots&=     \cdots\\
x_{n}/x_{1}&=(f_{n}/f_{1})(x/x_1)\;,
\end{cases}$$
where the division is in minplus sense, that is the subtraction.
Denoting $y=(x_{2}/x_{1},\cdots, x_{n}/x_{1})$ and
$g_{i-1}(y)=(f_{i}/f_{1})(0,y)$, the problem is reduced to  the
computation  of   the  fixed  point problem $y=g(y)$. This fixed
point gives the normalized eigenvector from which  the eigenvalue
is deduced by: $\lambda=f_{1}(0,y)$. We note that $g$ is a
non-homogeneous minplus function of $y$.

The fixed  point problem does not always have a solution, but
nevertheless, there are cases where we are able to find one:
\begin{itemize}
\item \emph{$f$ is affine in standard algebra}. In this case
$f(x)=Ax+b$. The homogeneity\footnote{$A\bar{1}=\bar{1}$ with
$\bar{1}$ the vector  with all  its  entries equal  to  1.}
implies that the kernel  of $A-I_{d}$  is not  empty.  When  this
kernel has  one dimension, the eigenvalue  is  equal to
$\lambda=pb$, where  $p$  is the  normalized ($p\bar{1}=1$) left
standard eigenvector of $A$ associated to the (standard)
eigenvalue $1$. But even in this case, all the standard
eigenvalues do not have  a module necessarily smaller than one,
and the dynamical system may be unstable. We note that when all
the entries of the matrix $A$ are nonnegative, $f$ is monotone
nondecreasing, but when there are positive entries and negative
entries, the system is not monotone. \item \emph{$f$ is minplus
linear}: $f(x)=A\otimes x$. In this case the system is monotone.
\item \emph{$f$ corresponds to stochastic control}. In this case
$f(x)=D\otimes(Hx)$ where $H$ is a standard matrix with rows that
define discrete probability laws. Such a matrix is called a
stochastic matrix. Then  the  dynamics $x^{k+1}=f(x^k)$ has the
interpretation of a dynamic programming equation associated to a
stochastic control problem, and  the eigenvalue is the optimal
average cost of the corresponding stochastic control problem. Note
that in this case, $x^k$ are components of all the dynamics which
can be written
$$x^{k+1}=A \boxtimes x^k, \text{ with }A=\begin{bmatrix} 0 & H \\
D & \varepsilon \end{bmatrix}\;.$$ \item \emph{$f$ corresponds to
stochastic games}. In this case $f=D_1\odot (D_2 \otimes (H x))$
where $\odot$ denotes the maxplus product (obtained by replacing
min by max in the minplus matrix product $\otimes$), and $H$ is a
stochastic matrix. This case corresponds to dynamic programming
equations associated to stochastic games. In this case, $f$ is
monotone.
\item \emph{$f$  has a particular  triangular structure} for example:
$$\begin{cases} x^{k+1}=A\otimes x^k\; , \\
y^{k+1}=B(x^k)\otimes y^k\; ,\end{cases} $$ with $B(x)$ additively
homogeneous of degree $0$. The system is not necessarily monotone. For such
systems, it is easy to find the eigenvalue and eigenvector by
applying the minplus algebra results. See~\cite{FARTh, Russ} for discussions and
generalizations.
\end{itemize}

In general, it is possible to compute the fixed point using the
Newton's method\footnote{ we have to  solve a piecewise linear
system of equations}. This method corresponds to the  policy
iteration when the dynamic programming  interpretation holds true.
For stochastic control problems, the policy iteration is globally
stable. In the game case, the policy iteration is only locally
stable.

We  may  have  unstable  fixed  points which  are  not  accessible
by integrating  the   dynamics. In  this case, the  eigenvalue, computable
by  the Newton's method, gives no information on the time asymptotic of
the system.  When all the  fixed points  are unstable,  we  may
have  a linear  growth of the state trajectories. This point is illustrated
by the chaotic tent dynamics example given in the next
section.

\subsection{Growth rate of homogeneous systems}
We define the \emph{growth rate}, $\chi(f)$, of a dynamic system
$x^{k+1}=f(x^k)$, where $x\in\Rbmin^n$, by the common limit
$\lim_k x_i^k/k$ of all the components $i$ when this limit exists.
In \cite{GG04} it has been proved, with a special definition of
connexity (satisfied for a system defined by $f(x)=A\boxtimes x$
as long as the graph $\GG(A)$ is strongly connected\footnote{Here
there is an edge from $i$ towards $j$ in $\GG(A)$ if $A_{ji}\neq
\varepsilon$ if it is a minplus edge or $A_{ji}\neq 0$ if it is a
standard one.}), that the growth rate and the eigenvalue of a
homogeneous and monotone system exist and are equal. Let us show
that chaos may appear and that the eigenvalue and growth rate which
exist are not equal on a system which is only homogeneous.

Let us consider the homogeneous dynamic system where $k$ is the time index:
$$\begin{cases}  x_{1}^{k+1}=x_{2}^k\;, \\
 x_{2}^{k+1}=(x_{2}^{k})^3/(x_{1}^k)^2\oplus 2(x_{1}^{k})^2/x_{2}^{k}\;.
\end{cases}$$

The corresponding eigenvalue problem is
$$\begin{cases} \lambda x_{1}=x_{2}\;, \\ \lambda x_{2}=x_{2}^3/ x_{1}^2\oplus 2x_{1}^2/x_{2}\;,
\end{cases}$$
where the minplus power exponent must not be confused with a time index.
\dessin{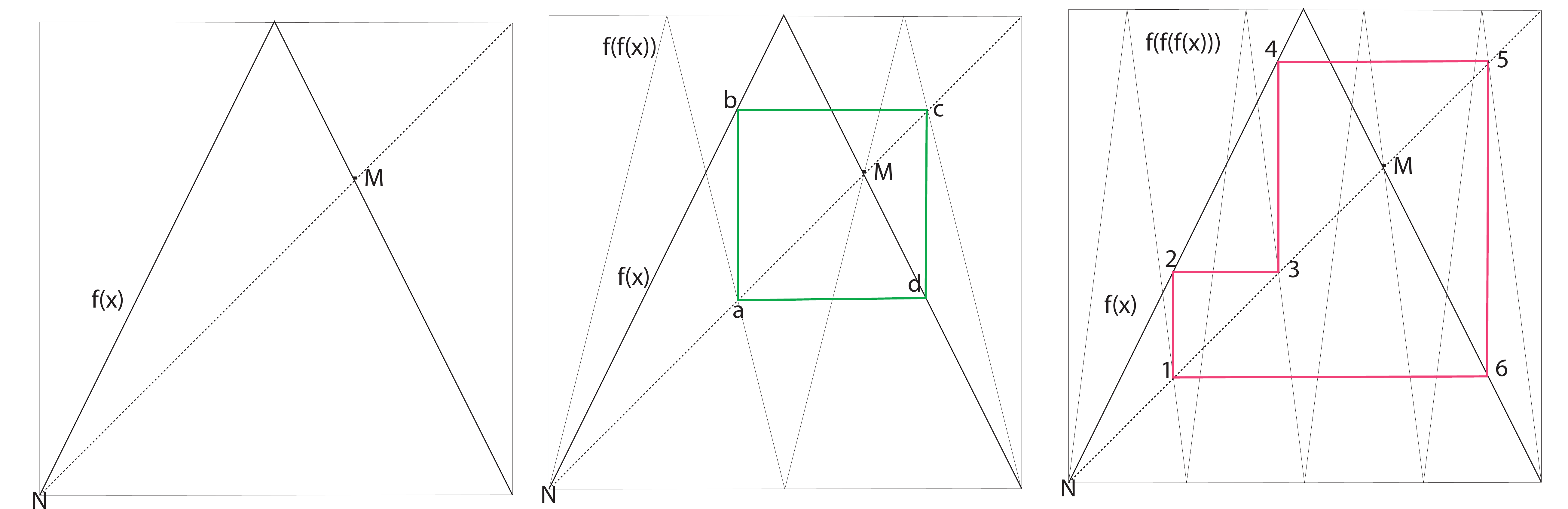}{0.30}{Cycles of the tent transformation. The abscissa $x_M$ of $M$ is a fixed point of $f$:
$x_M=f(x_M)$. The pair $(x_a,x_c)$ is a cycle of $f$ composed of two fixed points of $f\circ f$:
$x_c=f(x_a)$, $x_a=f(x_c)=f(f(x_a))$. The triplet $(x_1,x_3,x_5)$ is a circuit of $f$
composed of three fixed points of $f\circ f \circ f$: $x_3=f(x_1)$, $x_5=f(x_3)=f(f(x_1))$,
$x_1=f(x_5)=f(f(f(x_1)))$.}{tent}

The solution  is  $\lambda=y$  with $y=x_{2}/x_{1}$  satisfying the
equation  $$y=y^2\oplus 2/y^2,$$  which  has the  solutions $y=0$
and $y=\frac{2}{3}=0.66...$.  These  two  solutions  are  unstable
fixed points  of  the transformation  $f(y)=y^2\oplus  2/y^2$.
However, the system  $y^{k+1}= f(y^{k})$  is  chaotic since $f$ is
the tent transform (see~\cite{BERG} for a clear discussion of this
dynamics). In Figure-\ref{tent}, we show the graph  of $x\mapsto
f(x)$, $x\mapsto f(f(x))$ , $x\mapsto f(f(f(x)))$, their fixed
points, and   periodic trajectories.

In Figure-\ref{tentsimu}, we show a trajectory for an initial
condition chosen  randomly  with  the  uniform  law on  the  set
$\{(i-1)/10^5, i=1,\cdots,10^5\}$. The  diagonal line in the
picture is a decreasing order applied to the set
$\{y^k,\,k=1,\cdots,10^5\}$. This line shows that the invariant
empirical density is approximately uniform. In fact, with these
initial conditions, the trajectories are periodic with a possible
long period. \dessin{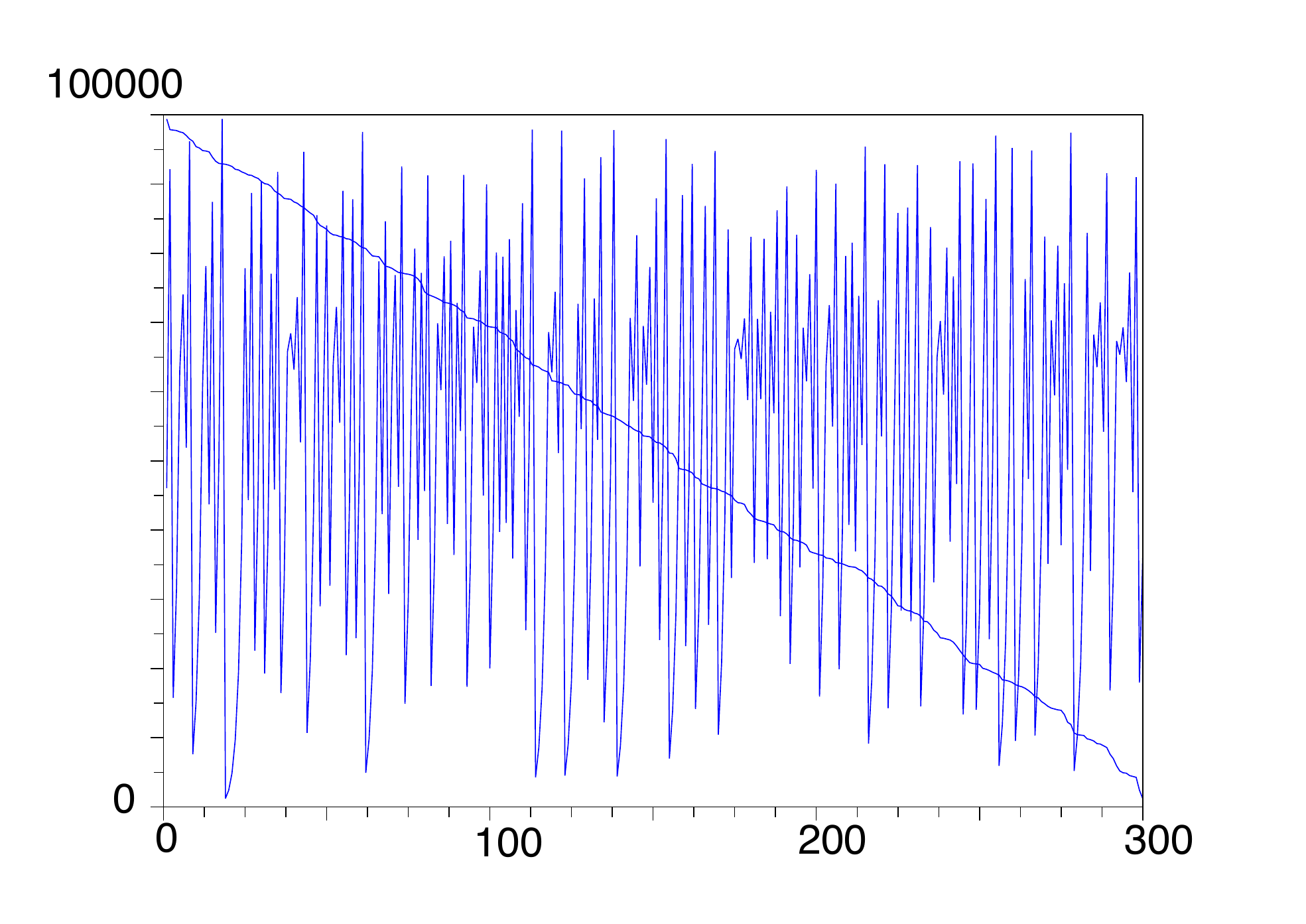}{0.45}{A tent iteration trajectory
($10^5y^k)_{k=1,.,300}$) and its arrangement in decreasing order
(almost diagonal line) .}{tentsimu} It has been proved that the tent
iteration has  a unique invariant measure absolutely  continuous
with respect to the Lebesgue measure (the uniform law on $[0,1]$).
Therefore, the system is ergodic.  The growth rate
$$(x_1^N-x_1^0)/N=\frac{1}{N}\sum_{k=1}^N(x_1^{k}-x_1^{k-1})\;,$$
can be computed by averaging, with respect to the uniform law, an
increase in one step: $f_1(x)-x_1$ with the standard notations
(that is $f_1(x)/x_1=y$ with the minplus notations). Therefore
$\chi(f)=\int_0^1ydy=0.5\;,$ for almost all initial conditions,
which is different from the eigenvalues ($0$ and $\frac{2}{3}$).

More generally, for a homogeneous system, we can write the system
dynamics:
$$x_1^{k+1}/x_1^k=f_1(x^k)/x_1^k=h(y^k),\;\; y^{k+1}=g(y^k),$$
with $y_{i-1}^{k}=x_i^k/x_1^k$ and  $g_{i-1}=f_i/f_1$ for
$i=2,\cdots,n$. As long as the $y^k$ belong to a bounded closed
(compact) set for all $k$, we remark (after Kryloff and
Bogoliuboff) that the set of measures:
$$\left\{P^N_{y^0}=\frac{1}{N}\left(\delta_{y^0}+\delta_{g(y_0)}+\cdots+\delta_{g^{N-1}(y^0)}\right),\;n\in\bbN\right\}\;,$$
(where $\delta_a$ denotes the Dirac mass on $a$) is tight.
Therefore, we can extract convergent subsequences which converge
toward invariant measures $Q_{y^0}$ that we will call
Kryloff-Bogoliuboff invariant measure.
Then we can apply the ergodic theorem at the sequence
$(y^k)_{k\in\bbN}$.  Application of this theorem shows that, for
almost all new initial conditions chosen randomly according to the
$Q_{y^0}$, we have~:
\begin{equation}\chi(f)=\lim_N\frac{1}{N}(x_1^{N}-x_1^0)=\lim_N \frac{1}{N}\left(\sum_{k=0}^{N-1}h(y^k)\right)=\int h(y)dQ_{y^0}(y)\;.\label{ergod}\end{equation}

It may happen that the initial
condition $y^0$ is transient; therefore, it is in the attractive
basin of $Q_{y^0}$ not in the support of $Q_{y^0}$. It would be
very useful to prove that $y^0$ is generic (in the sense of
Furstenberg\cite{Fur}), that is the limit exists for $y^0$. A
priori homogeneous systems do not have the uniform continuity
property required to prove the convergence of the Cesaro means for
all $y^0$.  The following classic example show the case of an ergodic 
system where with non generic points $f:x\in \bbT^1\rightarrow 2x\in \bbT^1$ with: $x^0=0.1\overset{2}{00}\overset{4}{1111}\overset{8}{00000000}\cdots$ where $\bbT^1$
the torus of dimension 1. In Figure-\ref{cesaro} we see that the Cesaro does not converge.
\dessin{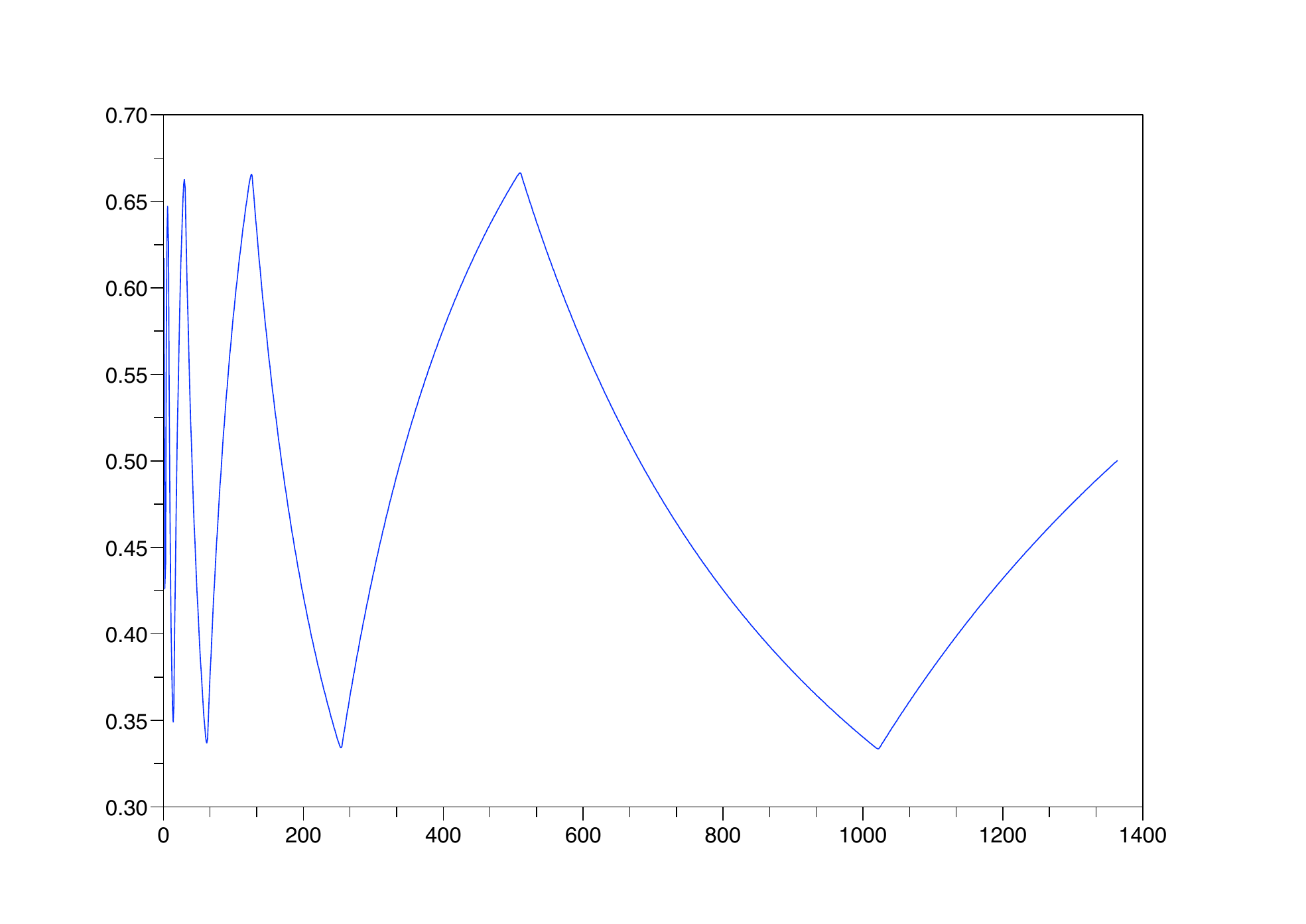}{0.4}{Plot of $S(n)$ with: $S(n)\triangleq \frac{1}{n}\sum_{k=0}^{n-1}x^k$.}{cesaro}

In the case where the compact set is finite, we can
apply the ergodicity results on Markov chains with a finite state
number to show the convergence of $P^N_{y^0}$ towards $Q_{y^0}$.
Instead of the subsequence convergence, this convergence proves
the genericity of $y^0$.

This discussion is summarized by the following result.
\begin{theorem}\label{GrowthRate}
For any additively 1-homogeneous dynamical systems $x^{k+1}=f(x^k)$, $x^0$ such that
$x^k_j/x^k_1$ stays bounded for all $j$ and $k$, there exists a measure on the initial
condition $x^0$ such that the growth rate exists for almost all the initial conditions.
\end{theorem}
We can also see \cite{AnKons} for construction of invariant
measures of stochastic recursions.

Coming back to the tent example, according to the  initial value
$y^0$, the  tent iterations $y^{k}$  stay  in  circuits  or follow
trajectories  without  circuit (possibly dense in $[0,1]$). For
example, assuming that the initial condition is
$y=\frac{2}{5}$, the trajectory is periodic of period 2. The
invariant measure is
$Q_{y^0}=\frac{1}{2}(\delta_{\frac{2}{5}}+\delta_{\frac{6}{5}})\;.$
The growth rate is $\frac{4}{5}$, which is again different from
the eigenvalues $0$ and $\frac{2}{3}$. Moreover, it can be shown
that for all initial conditions with a finite binary development
(this set contains all the float numbers of computers), the
trajectory stays in the unstable fixed point $0$ after a finite
number of steps. That is, for a dense set of initial conditions
the invariant measure is $\delta_0$ and the growth rate is $0$.

\section{Petri Net Dynamics}

\subsection{Autonomous Petri nets}
We present, in  min-plus-times  algebra,
timed continuous Petri nets with weights. The weights can be
negative and the numbers  of tokens  are not  necessary integer
(in  continuous Petri nets, what we call tokens are in fact fluid
amounts).

A  Petri net $\mathcal  N$ is  a graph  with two  sets of  nodes
(the \emph{transitions} $\mathcal Q$ (with $|\mathcal Q|$
elements) and the \emph{places}  $\mathcal P $  (with $|\mathcal
P|$ elements)) and two sorts of edges,  (the \emph{synchronization
edges} (from  a place  to a transition)  and the \emph{production
edges} (from  a transition  to a place)).

A  minplus  $|\mathcal  Q|\times  |\mathcal  P|$  matrix  $D$,
called \emph{synchronization\footnote{Decision matrix  in
stochastic control.} matrix} is  associated to  the
synchronization edges.   $D_{qp}=a_p$ if there exists an edge from
the place $p\in \mathcal P$ to the transition $q\in \mathcal Q$,
and  $D_{qp}=\varepsilon$ elsewhere, where $a_p$ is the
\emph{initial marking} of the place $p$, which is, graphically,
the number of tokens in $p$.  We suppose here that the sojourn
time in all the places  is one unit of time\footnote{When
different integer sojourn times are  considered, an equivalent
Petri net with a  unique sojourn time can be obtained by  adding
places and transitions and solving the implicit relations.}.

A  standard  algebra $|\mathcal  P|\times  |\mathcal  Q|$ matrix
$H$, called \emph{production\footnote{Hazard  matrix in stochastic
control} matrix}  is associated  with the  production edges.   It
is  defined by $H_{pq}=m_{pq}$  if there  exists  an edge  from
$q$ to  $p$, and  $0$ elsewhere, where $m_{pq}$ is the
\emph{multiplicity} of the edge \footnote{Here the multiplicity
appears only with the output transition edges. The multiplicity of
input transition edges is supposed to be always equal to one.
Looking at the more general case \cite{CGQ}, we  see that we do
not lose generality  by  doing so  (the dynamics class obtained is
the same)}.

Therefore, a Petri net is characterized by the quadruple:
$$(\mathcal P  ,\mathcal Q,H,D)\;.$$ A Petri net is a dynamic system  in which
the  token (fluid) evolution  is partially  defined by  the
transition firings, saying that a transition can  fire as soon as
all its upstream places contain a positive quantity  of tokens
(fluid) having stayed at least one unit of time.
When  a transition  fires, it  consumes a  quantity of  tokens (fluid)
equal  to the minimum  of all  the available  quantities being  in the
upstream  places.  Cumulating the  firings  done  up  to present  time
defines the \emph{cumulated transition  firing} of the transition. The
firing produces a quantity of  tokens (fluid) in each downstream place
equal to the  firing of the transition multiplied  by the multiplicity
of  the  corresponding  production  edge.  If  the  multiplicity  of  a
production edge, going from $q$ to  $p$, is negative, the firing of $q$
consumes tokens (fluid) of $p$.

A general Petri Net defines constraints on the transition firing.
Denoting by $q^{k}$
the cumulated firings  of transitions $q\in \mathcal{Q}$ up to instant $k$,  they satisfy the constraints:
\begin{equation}\min_{p\in
q^{in}}\left[a_p+\sum_{q\in p^{in}}m_{pq}q^{k-1}-\sum_{q\in p^{out}}q^{k}\right]=0,\;\forall q\in \mathcal{Q},\;
\forall k,\label{GPetri}
\end{equation}
where  $p\in \mathcal{P}$ is  a place of  the Petri Net;
$q^{in}\subset \mathcal{P}$  [resp.   $q^{out}\subset \mathcal{P}$]
denotes set of places upwards [resp. downwards]  the transition
$q$; and $p^{in}\subset \mathcal{Q}$ [resp. $p^{out}\subset
\mathcal{Q}$] denotes the set of transitions upwards [resp.
downwards] the place $p$.

Indeed, being at time $k$, we know (from the firing definition of
transitions) that after the firing (which is instantaneous), there
is at least one place upstream of any transition in which no token
entered before time $k-1$. For each transition $q$, the
equation~(\ref{GPetri}) computes the number of tokens which have
stayed at least one unit of time in each place $p\in q^{in}$. The
equation expresses that at least one place is empty and the others
have nonnegative numbers of tokens.

As long as there is more  than one edge leaving a place, the
trajectory of the system is not well-defined  because we do not
know the path of a token leaving this place.

In the case of a  \emph{deterministic} Petri net 
(generally called conflict free Petri net) where all the
places have only one  downstream edge, the dynamics are
well-defined, meaning there  is  no  token   consumption conflicts
between the transitions downstream of each place\footnote{In the
non-deterministic case, we have to specify the  rules  which
resolve  the  conflicts by, for example,  giving priorities to the
consuming transitions or by imposing ratios  to be followed. As
long as these  rules are added,   the  initial non-deterministic
Petri net becomes  a deterministic  one.}. Then, denoting by
$Q=(q^k)_{q\in \QQ, k \in \mathbb N}$ the vector of sequences of
cumulated firing quantities of transitions, and by $P=(p^k)_{p\in
\PP, k \in \mathbb N}$  the vector  of  sequences of  cumulated
token quantities arrived in the places at time $k$, we have:
  \begin{equation}\label{dynamic}
    \begin{bmatrix}P^{k+1} \\ Q^{k+1}\end{bmatrix}=\begin{bmatrix}0 & H \\ D & \varepsilon \end{bmatrix}
    \boxtimes\begin{bmatrix}P^{k+1}        \\       Q^{k}\end{bmatrix}
    \stackrel{\text{def}}{=}   \begin{bmatrix}H    Q^k   \\   D\otimes
    P^{k+1}\end{bmatrix}\; .
  \end{equation}
  
In Equation~(\ref{dynamic}), $P^{k+1}$ counts the number of tokens available (for firing in the
downstream transitions)
at time $k+1$ coming from the upstream transition firings.
It is obtained by summing the weighted firing numbers 
of the transitions upstream the places since the tokens are supposed to stay at least one unit of time in the places. 
 
The part about $Q^{k+1}$ in Equation~(\ref{dynamic}) tells that the transition firing numbers  at time $k+1$ are equal to the minimal token numbers in the places upstream
the transitions. To obtain these quantities we have only to add the token numbers in the place at initial time to the numbers entered
by the firings that are the entries of $P^{k+1}$.
 
From these dynamics, we deduce the dynamics of the cumulated firing
quantities by eliminating the place variables.  We deduce the
dynamics of the cumulated token quantities by eliminating the
transition variables.
$$Q^{k+1}=D\otimes (HQ^k),\quad P^{k+1}=H(D\otimes P^k)\; .$$

In  the case  of \emph{event  graphs} (particular  deterministic
Petri nets where  all the multiplicities $m_{pq}$  are equal to  $1$
and all the places  have exactly one edge  upstream), the dynamics
are linear in the minplus sense~:
$$Q^{k+1}=A \otimes  Q^k, $$ where  $A_{q'q}=a_p$ with $p$  the unique
place between $q$ and $q'$.

\subsection{Deterministic Petri nets}\label{DPetri}
By using negative weights and/or fixing a routing policy, it is possible to transform a Petri net with conflicts to a deterministic Petri net. 
Let us discuss
these points more precisely on the simple system given in the
first picture of Figure-\ref{homog}.
\dessin{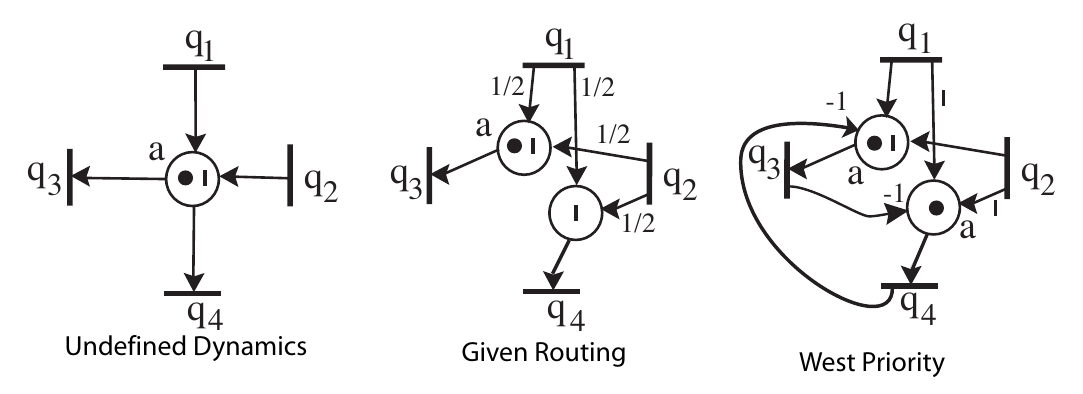}{1}{The nondeterministic Petri net, given in the
left figure, is made clear by: -- choosing a routing policy: 1/2
(to read in standard algebra) towards $q_3$, 1/2 towards $q_4$ in
the central figure, -- giving top priority to $q_3$ against $q_4$
in the right figure (where the time shift given by the sticks are
not anymore only in the places but also on the edges).}{homog}

The  incomplete dynamics  of this  system  can be  written in  minplus
algebra\footnote{Which means in standard algebra:
$q_{4}^{k}+q_{3}^k=a+q_{1}^{k-1}+q_{2}^{k-1}$.}:
\begin{equation}q^k_{4}q^k_{3}=aq^{k-1}_{1}q^{k-1}_{2}\;.\label{constr}\end{equation} Clearly $q_{3}$ and $q_{4}$
are not defined  uniquely. In the following two examples, we 
complete the dynamics in two different ways.
These two ways are equally useful for traffic
applications as we will see later.
\begin{itemize}
\item By specifying the routing policy (for example we choose arbitrarily that half of the total
tokens available are given to $q_3$ and half to $q_4$, see~\cite{CGQ} for results on the general routing case) \footnote{Which means in standard
algebra: $q^k_{4}=\frac{1}{2}\left(q_{1}^{k-1}+q_{2}^{k-1}\right), \quad q^k_3=a+q^k_4\;.$}:
$$ q^k_{4}=\sqrt{q^{k-1}_{1}q^{k-1}_{2}},\quad  q^k_{3}=aq^k_{4}\;.$$
The minplus product of the two equations gives the constraint~(\ref{constr}).
\item  By choosing  a priority  rule (top priority to $q_3$ against $q_4$)\footnote{Which means  in standard
algebra: $q^k_{3}=a+q^{k-1}_{1}+q^{k-1}_{2}-q^{k-1}_{4}, \quad q^k_{4}=a+q^{k-1}_{1}+q^{k-1}_{2}-q^k_{3}\;.$}:
$$ q^k_{3}=aq^{k-1}_{1}q^{k-1}_{2}/q^{k-1}_{4},\quad  q^k_{4}=aq^{k-1}_{1}q^{k-1}_{2}/q^k_{3}\;.$$
The last equation implies that the initial
constraint~(\ref{constr}) is satisfied. We see
that the negative weights on $q_4^{k-1}$ and on $q_3^k$ are
essential to express this priority. 
\end{itemize}
In the  two cases, we obtain a  \emph{degree one  homogeneous
minplus} system.

\subsection{Input-Output Petri nets}
We can define a Petri net with inputs and outputs in the following
way. We  partition  the transition  set  in three  parts
$(\VV,\QQ,\ZZ)$: the input set  $\VV$, the state set $\QQ$ and
the output set  $\ZZ$. We do  the same for the place  set
and we  get the three   parts  $(\UU,\PP,\YY)$.   The   inputs are
the  transitions [resp. places] without upstream edges. The
outputs are the ones without output edges.  Then the dynamics can
be rewritten:
\begin{equation}\begin{bmatrix}
           P^{k+1}\\ Q^{k+1}\\ Y^{k+1}\\ Z^{k+1}
   \end{bmatrix}=
          \begin{bmatrix}
           0 & A & 0 & B\\ C & \varepsilon & D & \varepsilon\\ 0 & E &
           0 & 0\\ F & \varepsilon & \varepsilon & \varepsilon
           \end{bmatrix}\boxtimes
           \begin{bmatrix}
           P^{k+1}\\ Q^k\\ U^{k+1}\\ V^k
          \end{bmatrix}
          \triangleq
          \begin{bmatrix}
               A  Q^k  +   B  V^k\\  C\otimes  P^{k+1}\oplus  D\otimes
           U^{k+1}\\ E Q^k \\ F\otimes P^{k+1}
           \end{bmatrix}\; . \label{PetriProd} \end{equation}
These dynamics, denoted by $S$ and defined by the matrices
$(A,B,C,D,E,F)$, associate the output signals
$(Y^k,Z^k)_{k\in\mathbb N}$ to the input signals
$(U^k,V^k)_{k\in\mathbb N}$. We write $(Y,Z)=S(U,V)\;$. These
systems can be considered as a special case of an extension of the
generalized system described in the minplus section.  In this
special case, some blocks are null and an implicit part can
appear.

The three operations (parallel, series, and feedback compositions)
described previously can be extended easily to this case.
\section{Traffic Application}

\subsection{Traffic on a circular road}
Let us recall the simplest model to derive the fundamental traffic
diagram on a single  road. The best way to obtain this diagram is to  study the
stationary regime on  a circular  road  with  a given  number  of
vehicles\footnote{We consider that  the stationary regime  of
the circular road  is reached locally  on  a standard  road  when
its density  is  constant at  the considered  zone.}. We present
two ways  to obtain  this diagram: one is by logical  deduction
from an exclusion  process  point of view 
(it shows clearly the presence of two distinct phases), and the other by
computing the eigenvalue of a minplus system  derived from  a
simple Petri net modeling  of the road (this way will be extended to the case of roads with junctions). In the following, the road
is cut in $m$ sections which can    contain   at   most one
vehicle. \dessin{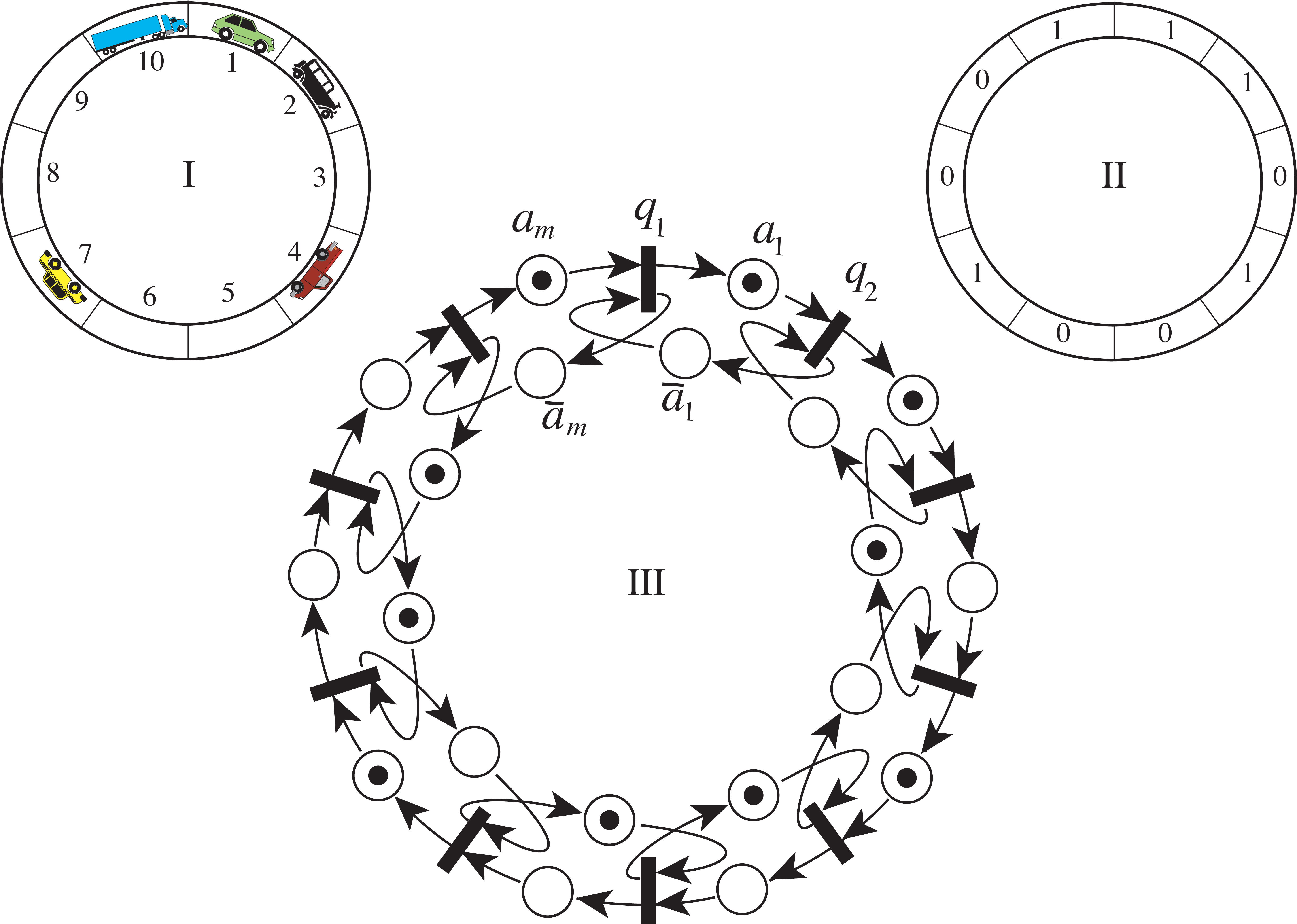}{0.2}{On the top-left side, a circular road
cut in sections. On the top-right side, its exclusion process,
where 1 means that the section is occupied. On the middle, its
Petri net representation where the ticks (1 time delay) in each
place are not represented.}{ruta}
\subsection{Exclusion process modeling}
Following~\cite{BLA}, we can consider the dynamical system defined
by the rule $10\rightarrow  01$ applied  to a binary  word $w$.
The word $w_k$ describes  the vehicle positions at instant $k$  on
a road cut in  sections (each  bit representing a section, 1
meaning occupied and 0 meaning free, see II in Figure-\ref{ruta}).
Let us take an example:
\begin{align*}
w_{1}  &=1101001001,\quad  w_{2}=1010100101,  \\  w_{3}  &=0101010011,
\quad w_{4}=1010101010, \\ w_{5} &=0101010101. \\
\end{align*}

Let us define: -- the \emph{density} $\rho$ by the number of
vehicles $n$ divided by the number of sections $m$: $\rho=n/m$,
-- the flow $\varphi(t)$  at time $t$  by the  number of  vehicles
going  one step forward  at  time  $t$  divided  by  the number of
sections.   Then  the \emph{fundamental traffic diagram} gives the
relation between $\varphi(t)$ and $\rho$.

If $\rho \leq  1/2$, after a transient period,  all the vehicle groups
split off,  and then all the  vehicles can move  forward without other
vehicles in the way, and we have:
  $$\varphi(t)=\varphi=n/m=\rho\;.$$  If  $\rho  \geq 1/2$,  the  free
place groups split  off after a finite time  and move backward without
other free places in the way.  Then $m-n$ vehicles move forward and we
have:
$$\varphi(t)=\varphi=(m-n)/m=1-\rho\;.$$
\begin{theorem}
$$\exists  T:\;\forall  t\geq T\quad  \varphi(t)=\varphi=\begin{cases}
\rho &  \mathrm{if} \;  \;\rho\leq 1/2\;, \\  1-\rho &  \mathrm{if} \;
\;\rho\geq 1/2\; .
\end{cases}$$
\end{theorem}

\begin{proof}
This result has been proved in \cite{BLA}. Let us give the idea of
the proof. We only have to prove that after a finite number of
steps, the separations of vehicles or holes appear.

Assuming that the density is not greater than 1/2, let us look
what happens at a cluster of at least two vehicles denoted by A.
There are two cases:
\begin{itemize}
\item The cluster behind A is separated from A by one hole, we
have the configuration 10A0. Then at the next step we have 0A01.
The cluster A has gone backward by one place. \item The cluster
behind is separated from A by more than one hole, we have the
configuration 00A0. Then at the next step, the size of A has been
reduced by one.
\end{itemize}
Then after a finite number of steps, bounded by the number of
places, the size of each cluster of vehicles is reduced strictly.
Indeed, the density being not greater than 1/2, and individual
vehicles going forward, there is cluster of vehicles that must
meet a cluster of at least two holes.

For a density larger than 1/2, we follow the cluster of holes
instead of vehicles. We can show by the same arguments that their
sizes decrease or they go forward. Therefore after a finite number
of steps the holes are separated by vehicles.
\end{proof}

\dessin{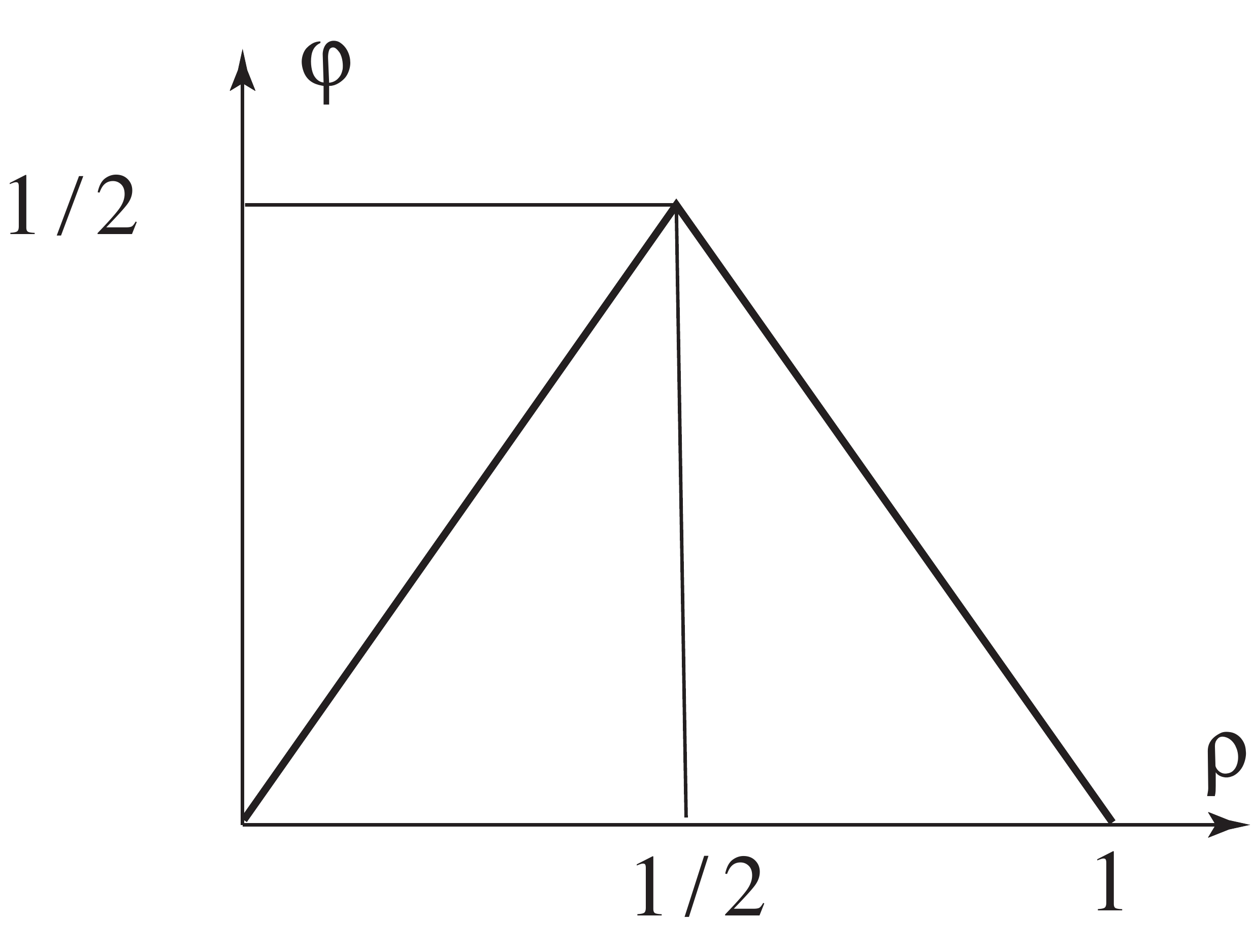}{0.2}{The fundamental traffic diagram showing the
dependence of the average car flow on the car density.}{fundamental}
\subsection{Event Graph modeling}
Consider  the  Petri  net  given  in III  of  Figure-\ref{ruta}
which describes the  same dynamics in  a different way. In  fact,
this Petri net is an event graph. Therefore, the dynamics are linear
in minplus algebra.  The number of vehicles entered  in the section $s$
before time $k$ is  denoted $q_{s}^k$.  The  initial vehicle
position is given by the booleans $a_{s}$ which take the value $1$
when the cell contains a vehicle and $0$ otherwise.

We use the notation $\bar{a}=1-a$. Then, the dynamics are given by:
$$q_{s}^{k+1}=\min\{a_{s-1}+q_{s-1}^k,\bar{a}_{s}+q_{s+1}^k\}\;,$$
which can be written linearly in minplus algebra:
$$q_{s}^{k+1}=a_{s-1}q_{s-1}^k\oplus\bar{a}_{s}q_{s+1}^k\;,$$
where the index addition is done modulo $m$.

%
\begin{theorem}
The average transition  speed (car flow) $\varphi$ depends on the
car density $\rho$ according to the law:
\[\varphi=\min(\rho,1-\rho)\; \square\]
\end{theorem}
\begin{proof}
This  event graph has  three kinds of elementary  circuits:
the outside circuit  with average weight $n/m$; the inside circuit
with average  weight $(m-n)/m$; and the circuits on which some
steps are made forward and then back with the average mean 1/2.
Therefore, using Theorem~\ref{eigenth}, its eigenvalue is
\[\varphi=\min(n/m,(m-n)/m,1/2)=\min(\rho,1-\rho)\; ,\] which  gives
the average  speed as a  function of the  vehicle density since the
minplus eigenvalue is equal to $\lim_{t}x_{i}(t)/t=\varphi$ for all
$i$. \end{proof}

\subsection{Traffic on two roads with one junction}
Before discussing the dynamics of the town, let us study in detail
the case of two circular roads with a junction (see the top-right
side of Figure-\ref{junctionsimp}).

A   first   trial   is   to   consider  the   Petri   net   given
in the middle of Figure-\ref{junctionsimp}.   This Petri net is
not an event graph.   It    is   a    general   non deterministic
Petri   Net. \dessin{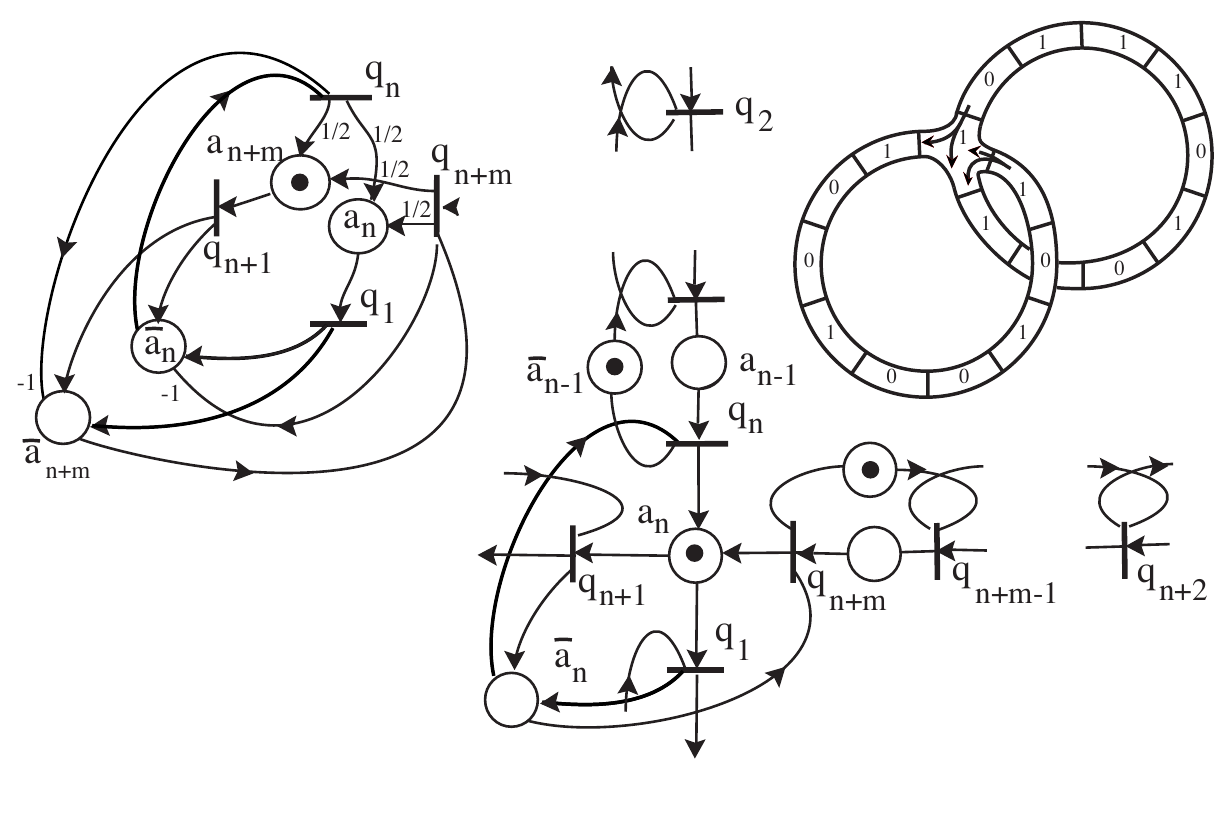}{0.9}{A junction with two
circular roads cut in sections (top-right), its Petri net
simplified modeling (middle) and the precise modeling of the
junction (top left). The ticks representing the time delays
present in each place are not represented.}{junctionsimp}

We   can    write   the   dynamics    of   this   Petri    net
using Equation-\ref{GPetri},  but these equations  do not uniquely
determine the trajectories of  the system. We have two  places
with two outgoing edges, $a_n$ and $\bar{a}_n$.  At place $a_n$,
we have to specify the \emph{routing policy} giving the proportion
of cars going towards $y_2$  and the proportion going  towards
$y_3$. At place $\bar{a}_n$, we  may follow \emph{the first
arrived the first served rule with the right priority} when two
vehicles  want to enter in the junction simultaneously. Using
Petri net  with negative weights, we obtain the Petri Net of
Figure-\ref{junctionsimp} where the junction is described
precisely in the top-left part of the figure. In the place $a_n$ [resp $a_{n+m}$] are accumulated
the cars  going towards West [resp. towards South]. In the place $\bar{a}_n$
[resp. $\bar{a}_{n+m}$] are accumulated the authorizations to enter in the junction from North
[resp. from East]. These authorizations are obtained by subtracting the authorizations to enter from East
[resp. from North] from the total number of authorizations to enter in the junction. See  Section~\ref{DPetri}
to have a precise description of the management of the priority rule. We remark that we have obtained a deterministic
(conflict free) Petri Net with some negative weights. 


The dynamics describing the evolution of the vehicle number entered in section
$s$ before time $k$ denoted $q_s^k$ can be obtained immediately from the Petri net
of the top left corner of Figure-\ref{junctionsimp} using the minplus notation discussed in Section-\ref{MaxPlus}.

\begin{equation}\begin{cases}
q_{i}^{k+1}=a_{i-1}q_{i-1}^k\oplus\bar{a}_iq_{i+1}^k,\;i\not=1,n,n+1,n+m,
\\                   q_{n}^{k+1}=\bar{a}_nq_1^kq_{n+1}^k/q_{n+m}^k\oplus
a_{n-1}q_{n-1}^k\;,\\
q_{n+m}^{k+1}=\bar{a}_{n+m}q_1^kq_{n+1}^k/q_{n}^{k+1}\oplus
a_{n+m-1}q_{n+m-1}^k\;,                                                \\
q_{1}^{k+1}=a_{n}\sqrt{q_{n}^kq_{n+m}^k}\oplus\bar{a}_1q_2^k\;,\\
q_{n+1}^{k+1}=                       a_{n+m}\sqrt{q_{n}^kq_{n+m}^k}\oplus
\bar{a}_{n+1}q_{n+2}^k\;,\\
\end{cases}\label{2Ddyna}\end{equation}
where the entries satisfy the following constraints (written with the standard notations):
\begin{itemize}
\item $0\leq a_i\leq 1$ for $i=1,\cdots,n+m$. These initial
markings give the presence, $1$, or absence, $0$, of a vehicle in
the road sections. However, here we see vehicles as fluid and can
relax this integer constraint. Moreover, the transition firing
cuts the tokens, and if we want to see the systems after some
firings, there is not necessarily an integer number of tokens in a
place. Therefore, it is better to accept real numbers of tokens
belonging to the $[0,1]$ interval; \item $\bar{a}_i=1-a_i$ for
$i\neq n,n+m$ they give the initial free spaces in the places;
\item $a_n+a_{n+m}\leq 1$ the maximum number of cars in the
junction is 1; \item $\bar{a}_n=\bar{a}_{n+m}=1-a_n-a_{n+m}$ give
the free place in the junction.
\end{itemize}

We remark that the system is homogeneous of degree 1 and that it is easy to write
the dynamics using the generalized matrix product. Indeed, all the ``monomials'' appearing in the right hand side of System~\ref{2Ddyna} (for example $\sqrt{q_{n}^kq_{n+m}^k}$) are linear in the standard algebra and can be computed by a standard matrix product from the $q$ vector. Then it is easy to obtain the complete right hand side from all the appearing ``monomials''
by a minplus matrix product. 
This matrix form is very useful
to simulate the system in ScicosLab~\cite{ScicosLab} since the minplus matrix product is implemented in this software. We remark also that the ``monomials'' like $q_1^k q_{n+1}^k/q_{n+m}^k$ introduce
negative entries in the standard matrices. Therefore the dynamics is not monotonic.


We can prove the existence of a growth rate thanks to the following two
results.

\begin{theorem}
  The   trajectories  of   the   states  of   the  junction   dynamics
  (\ref{2Ddyna}),  $(q^k_{i})_{k\in  \bbN}$,  starting  from  $0$,  are
  nonnegative and nondecreasing for all $i\; \square$
\end{theorem}
\begin{proof}
  Computing   $q^1$  using  the   fact  that   all  the   $a_{i}$  and
  $\bar{a}_{i}$ are  nonnegative, it  is clear that  $q^1\geq 0$.  Let us
  prove by induction that the trajectories are nondecreasing. It is true
  for $k=1$.  We suppose  that it  is true for  $k$, that  is $q^k\geq
  q^{k-1}$), and we prove  that it is also true for $k+1$, that is $q^{k+1}\geq q^k$. Rewrite~(\ref{2Ddyna})
  $q_i^{k+1}=f_i(q^k)$ for $i=1,\cdots,n+m$. The functions
  $f_i$ for $i\not=n,n+m$ are nondecreasing. Therefore, for such an $i$,
  we have:
  $$q^{k+1}_i=f_i(q^k)\geq f_i(q^{k-1})\geq q^k_i\;,$$ using first the induction
  hypothesis and then the dynamics definition.

  Let us prove the nondecreasing property of $q_n$.
  \begin{itemize}
  \item
  If $q_{n}^{k+1}= a_{n-1}q_{n-1}^k$ we have
  $$q^{k+1}_n=a_{n-1}q_{n-1}^k\geq a_{n-1}q_{n-1}^{k-1}\geq f_n(q^{k-1})= q^k_n\;.$$
  \item
  If  $q_{n}^{k+1}=\bar{a}_n q_1^k q_{n+1}^k/q_{n+m}^k$, using the dynamics, we have $q_{n+m}^k\leq \bar{a}_{n+m} q_1^{k-1}
     q_{n+1}^{k-1}/q_{n}^{k}$. Therefore, $$q_{n}^{k+1}\geq    \bar{a}_n q_{n}^kq_1^k q_{n+1}^k/\bar{a}_{n+m}
     q_1^{k-1}  q_{n+1}^{k-1}$$ which  gives $q_{n}^{k+1}\geq  q_{n}^k$
     using   the   induction  hypothesis and the assumption
     $\bar{a}_n =\bar{a}_{n+m}$.
 \end{itemize}

 The nondecreasing property of $q_{n+m}$ is proved in the same way.
  If $q_{n+m}^{k+1}= a_{n+m-1}q_{n+m-1}^k$, we have
  $q^{k+1}_{n+m}=a_{n+m-1}q_{n+m-1}^k\geq a_{n+m-1}q_{n+m-1}^{k-1}\geq f_{n+m}(q^{k-1})= q^k_{n+m}\;.$
  If $q_{n+m}^{k+1}=\bar{a}_{n+m} q_1^k q_{n+1}^k/q_{n}^{k+1}$,  we have $q_{n}^{k+1}\leq \bar{a}_{n} q_1^{k}
  q_{n+1}^{k}/q_{n+m}^{k}$ using the
  dynamics.  Therefore,
  $q_{n+m}^{k+1}\geq    \bar{a}_{n+m} q_{n+m}^{k}/\bar{a}_{n}$,
  which gives the result using $\bar{a}_n =\bar{a}_{n+m}$.
\end{proof}

\begin{theorem}\label{Bounded}
  The distances between any pair of states stay bounded.
  $$\exists c_{1}: \sup_{k}|q^k_{i}-q^k_{j}|\leq c_{1}, \forall i,j.$$
  Moreover
  $$\forall   T,   \exists  c_{2}:   \sup_{k}|q^{k+T}_{i}-q^k_{i}|\leq
  c_{2}T, \forall i\; \square$$
\end{theorem}
\begin{proof}
  This result comes from  the following inequalities (written in minplus algebra) deduced from the
  dynamics and the nondecreasing property of the trajectories:
  $$q_{i}^{k+1}\leq a_{i-1}q_{i-1}^{k}\;,i\neq 1,n+1$$
  $$q_{1}^{k+1}\leq  a_n\sqrt{q_{n}^k  q_{n+m}^k}\leq  a_n\sqrt{q_{n}^{k+n}q_{n+m}^k}\leq
  a_n\sqrt{b_1^{n-1}q_{1}^{k+1}q_{n+m}^k}\Rightarrow   q_{1}^{k+1}\leq
  a_{n}b_1^nq_{n+m}^k\;,$$
  with $b_j^k=\bigotimes_{i=j}^k a_{i}$.
  $$q_{n+1}^{k+1}\leq a_{n+m}\sqrt{q_{n}^k q_{n+m}^k}\leq a_{n+m}\sqrt{q_{n}^k q_{n+m}^{k+m}}\leq
  a_{n+m}\sqrt{q_{n}^kb_{n}^{n+m-1}q_{n+1}^{k+1}   }   \Rightarrow   q_{n+1}^{k+1}\leq
  a_{n+m}b_{n}^{n+m}q^k_{n}\;.$$ Therefore we have:
  $$q_n^k\leq  a_{n-1}q_{n-1}^{k-1}\leq\cdots \leq (b_1^{n})^2q_{n+m}^{k-n}\leq
  (b_1^n)^2a_{n+m-1}q_{n+m-1}^{k-n-1}\leq\cdots\leq  (b_1^{n+m})^2q_n^{k-n-m}.$$  The result
  follows from these inequalities which give bounds for all the the distances between two states
  and between the same state but at different times.
\end{proof}

By using Theorem~\ref{Bounded} and Theorem~\ref{GrowthRate}, we give an existence theorem
of the growth rate of the system~(\ref{2Ddyna}). The growth rate has the interpretation of
the average traffic flow.
\begin{theorem} There exists an initial distribution on $(q_j^0/q_1^0)_{j=2,n+m}$ the Kryloff-Bogoljuboff invariant measure   such that the average flow $$\chi=\lim_{k}q^k_{i}/k,  \;\;  \forall  i  \;  ,$$
of the dynamical system (\ref{2Ddyna}) exists almost everywhere $\square$
\end{theorem}

This result is not completely satisfactory. We would like to have
the existence of the growth rate for the initial condition
$q_i^0=0$ for all $i$. Nevertheless, we note that the numerical
approximation obtained by simulation of the growth rate always
exists. Indeed, in this case, the number of the approximated
states obtained by floating number approximation is finite since
they belong to a bounded set\footnote{We suppose here that the
approximation does not destroy the fact that the states stay in a
bounded set .}. The existence of the growth rate can be proved
using the Cesaro-convergence of the
probability to be in a state toward an invariant measure in the case of finite Markov chains. 
At this point, it is also useful to recall that the continuous Petri net
model used here is an approximation of a more realistic discrete Petri net.

\begin{theorem}
  Starting from $q^0=0$, if an average growth rate $\chi$ exists, then $\chi\leq 1/4$.
\end{theorem}

\proof From the dynamics~(\ref{2Ddyna}) we have~:
\begin{align*}
  & q_{n+m}^{k+1}\leq \bar{a}_{n+m} q_1^k q_{n+1}^k/q_n^{k+1},\\
  & q_1^{k+1}\leq a_n\sqrt{q_n^k q_{n+m}^k},\\
  & q_{n+1}^{k+1}\leq a_{n+m}\sqrt{q_n^k q_{n+m}^k}.
\end{align*}
then by summing these inequalities we get~:
$$(q_1^{k+1}-q_1^k)+(q_n^{k+1}-q_n^k)+(q_{n+1}^{k+1}-q_{n+1}^k)+(q_{n+m}^{k+1}-q_{n+m}^k)\leq 1.$$
Hence, when $k\to +\infty$, and taking into account Theorem-\ref{Bounded}, we obtain $4\chi \leq 1$ \endproof

\subsection{Eigenvalue Existence of the Junction Dynamics}

   Let us  consider the eigenvalue problem associated to the  dynamics  (\ref{2Ddyna}).
   It is defined as finding $\lambda$ and $q$ such that:
   \begin{equation}\begin{cases}
\lambda q_{i}=a_{i-1}q_{i-1}\oplus\bar{a}_iq_{i+1},\;i\not=1,n,n+1,n+m,
\\                  \lambda q_{n}=\bar{a}_n q_1 q_{n+1}/q_{n+m}\oplus
a_{n-1}q_{n-1}\;,\\
\lambda q_{n+m}=\bar{a}_{n+m}q_1q_{n+1}/(\lambda q_{n})\oplus
a_{n+m-1}q_{n+m-1}\;,                                                \\
\lambda q_{1}=a_{n}\sqrt{q_{n}q_{n+m}}\oplus\bar{a}_1q_2\;,\\
\lambda q_{n+1}=                       a_{n+m}\sqrt{q_{n}q_{n+m}}\oplus
\bar{a}_{n+1}q_{n+2}\;,\\
\end{cases}\label{Eigel}\end{equation}
with: $0\leq a_i\leq 1$ for $i=1,\cdots,n+m$, $\bar{a}_i=1-a_i$ for $i\neq n,n+m$,
$a_n+a_{n+m}\leq 1$ and $\bar{a}_n=\bar{a}_{n+m}=1-a_n-a_{n+m}$.

The eigenvalue problem can be solved explicitly.
\begin{theorem} \label{maintheo} The nonnegative eigenvalues $\lambda$ as a function of the density $d$,
written in the standard algebra, is given by:
$$\begin{array}{|c|c|c|c|c|}
\hline
d  & 0\leq d \leq \alpha & \alpha \leq d \leq \beta
&\min(\beta,\gamma) < d < \max(\beta,\gamma) & \gamma\leq d \leq 1 \\
 \hline
 \lambda & (1-\rho)d & 1/4 & (r-(1-\rho)d)/(2r-1+2\rho)& 0\\
 \hline
\end{array}$$
 with $N=n+m$, $\rho=1/N$, $r=m/N$ and $d=\big(\sum_{i=1}^{n+m}a_i\big)/(N-1)$ the density of vehicles,
 $\alpha=1/4(1-\rho)$,  $\beta=(r+1/2-\rho)/2(1-\rho)$ and $\gamma=r/(1-\rho)\; \square$
\end{theorem}

If $r>1/2$ when $N$ is large, the positive eigenvalue is unique and
has the simple approximation:
$$\lambda\simeq\max\left\{\min\left\{d,\frac{1}{4},
               \frac{r-d}{2r-1}\right\},\; 0 \right\}.$$

\begin{proof}
The whole proof is available in~\cite{ARX} (where the role of $m$ and $n$ are inverted). 
Only a sketch of the proof is given in this section.  The proof has two parts. 
The first part consists of
reducing the problem to a generalized eigenvalue problem  in a four
dimensional space (this part is easy and is not given here, it is available
in~\cite{ARX}). The second part consists of a verification of the
generalized minplus eigenvalue system of equations since we 
give explicit formulas for all the eigenelements (this part is
given in the Appendix, and is also available in details
in~\cite{ARX}). The eigenelements have been
obtained by solving explicitly the homogeneous affine systems with
five unknowns achieving the minimum in the reduced system. 
Numerical simulations have suggested the affine system achieving the minimum. 
Knowing it, we had only to verify the inequalities proving that it achieves
actually the minimum.

Since we want the eigenvalue as a function
of the density, we cannot use a numerical approach. Instead, we have
to find explicit formulas which are piecewise affine. This is the main
difficulty of the problem.

After verifying that $\lambda\leq1/4$, by elimination of
$q_i,\;i\neq 1,n,n+1,n+m,$ thanks to the minplus linearity of the
first equation of~(\ref{Eigel}), we obtain the
closed set of equations defining $q_i$, $i=1,n,n+1,n+m$; see Lemma~\ref{Reduced} in the appendix:
\begin{equation}  \begin{cases}
    q_n=(\bar{a}_n\mpdiv\lambda) q_1 q_{n+1}\mpdiv q_{n+m} \oplus (b_n\mpdiv\lambda^{n-1}) q_1\;,\\
    q_{n+m}=(\bar{a}_{n+m}\mpdiv\lambda^2) q_1 q_{n+1}\mpdiv q_n \oplus (b_m\mpdiv\lambda^{m-1}) q_{n+1}\;,\\
    q_1=(a_n\mpdiv\lambda) \sqrt{q_n q_{n+m}} \oplus (\bar{b}_n\mpdiv\lambda^{n-1}) q_n\;,\\
    q_{n+1}=(a_{n+m}\mpdiv\lambda) \sqrt{q_n q_{n+m}} \oplus (\bar{b}_m\mpdiv\lambda^{m-1}) q_{n+m}\;,
  \end{cases}\label{RedEig} \end{equation}
   where $b_n=\bigotimes_{i=1}^{n-1}a_i$ is the number of cars in the street with priority;
   $\bar{b}_n=\bigotimes_{i=1}^{n-1}\bar{a}_i$ is the number of free places in the street with
   priority; $b_m=\bigotimes_{i=n+1}^{n+m-1}a_i$ is the number of cars in the street without priority;
   and $\bar{b}_m=\bigotimes_{i=n+1}^{n+m-1}\bar{a}_i$ is the number of free places in the street without priority.

Observing from numerical simulations that four phases exist, it is
possible to specify their domains analytically, and give their
physical interpretations by observing the asymptotic regimes:
\begin{itemize}
\item \emph{Free moving}. When the density is small, $0 \leq d
\leq \alpha$, after a finite time, all the vehicles move freely.
\item \emph{Saturation}. When $\alpha \leq d \leq \beta$, the
junction is used at its maximal capacity without being bothered by
downstream vehicles. \item \emph{Recession}. When $\beta< d <
\gamma$, the crossing is fully occupied, but vehicles sometimes
cannot leave the crossing because the road they want to enter is
crowded. When $\gamma<\beta$, three eigenvalues exist on the
interval $[\gamma,\beta]$. In this case, the system is blocked.
\item \emph{Freeze}. When $\gamma\leq d\leq 1$, the road without
priority is full of vehicles. No vehicle can leave this road while
the vehicle being in the junction wants to enter.
\end{itemize}

Using these understanding of the phases we have been able
to guess where the minimum are achieved and then to find the solution 
by solving a standard linear system of 4 equations and 5 unknowns for each phase.
In the Appendix, we show that Table~\ref{Array} (written in standard
algebra) gives the eigenvalue and eigenvector formulas as a function
of the density $d$ for each phase.
\begin{table}[htb]
$$\begin{array}{|c||c|c|c|c|}
\hline
  & 0\leq d \leq \alpha & \alpha \leq d \leq \beta
&\min(\beta,\gamma) < d < \max(\beta,\gamma) & \gamma\leq d \leq 1 \\
 \hline \hline
 \lambda & (1-\rho)d & \frac{1}{4} & \frac{r-(1-\rho)d}{2r-1+2\rho}& 0\\
 \hline
 q_n& b_n-(n-1)\lambda  & b_n-(n-1)\lambda  & b_n -(n-1)\lambda & \bar{a}_{n+m}+\bar{b}_m \\
 \hline
 q_{n+m} & (n+1)\lambda-2a_n-b_n & (n+1)\lambda-2a_n-b_n &(n+1)\lambda-2a_n-b_n & -2a_n-\bar{a}_{n+m}-\bar{b}_m \\
 \hline
 q_1& 0  & 0  & 0 & 0 \\
 \hline
 q_{n+1}& a_{n+m}-a_n  & a_{n+m}-a_{n} &4\lambda-1+a_{n+m}-a_n & -2a_n-\bar{a}_{n+m}\\
 \hline
\end{array}$$
\caption{The eigenvalues and eigenvectors as function of the car densities.}
\label{Array}
\end{table}

\end{proof}

In Figure-\ref{LoiFond}, we show the fundamental diagram obtained by
simulation (using the maxplus arithmetic of the ScicosLab software
\cite{ScicosLab}) for a particular relative size $r$ of the  two
roads and the eigenvalue $\lambda$ given in Table~\ref{Array}. We
see clearly the four phases described in the proof of the eigenvalue
formula. On this figure, we see also  that the growth rate and the eigenvalue 
are very close to each other at least for three phases (among four phases).
\dessin{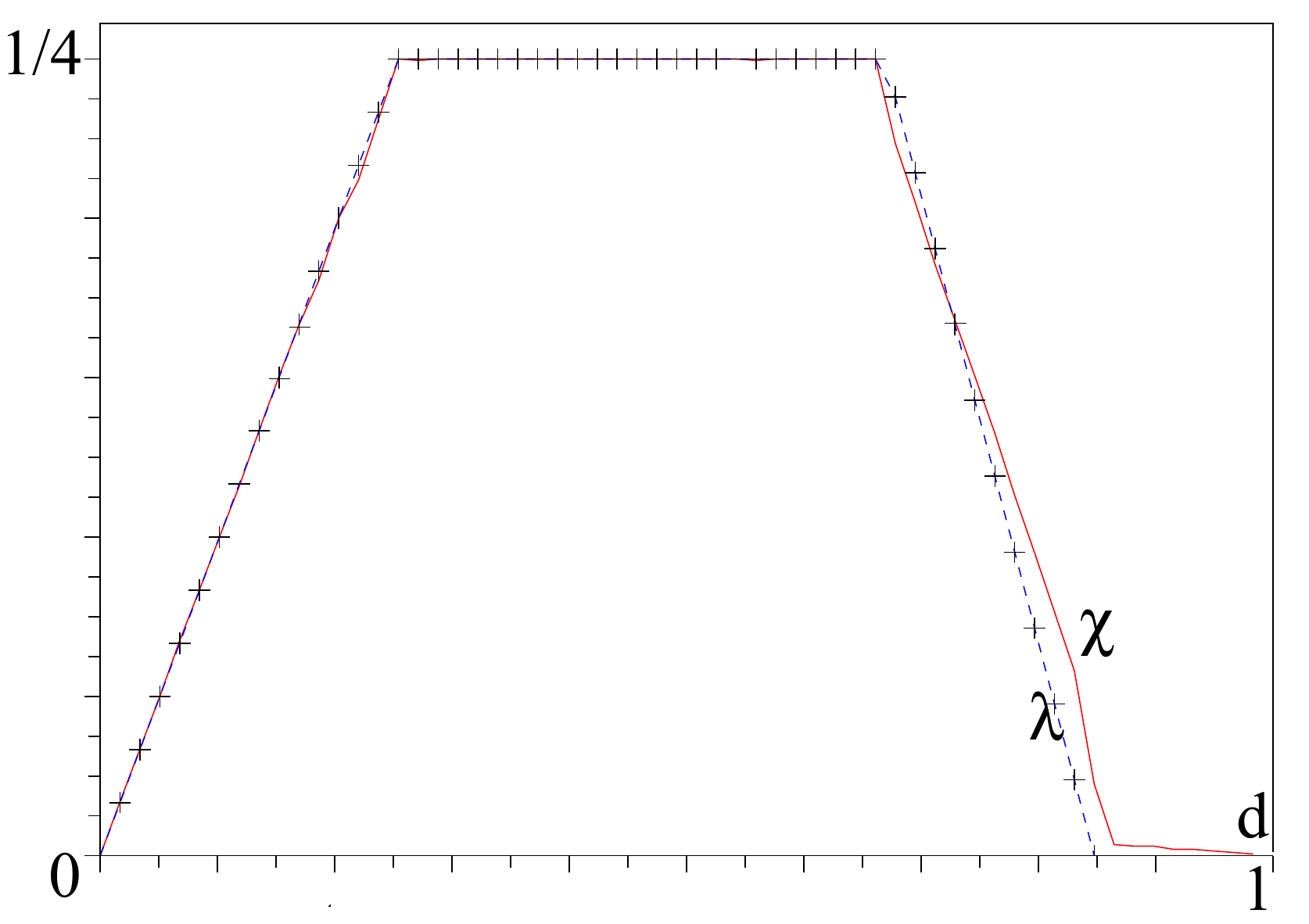}{0.40}{The traffic fundamental diagram
$\chi(d)$ when $r=5/6$ (continuous line) obtained by simulation
and its comparison with the eigenvalue $\lambda(d)$ given in
Table~\ref{Array}.}{LoiFond}
\subsection{Regular City Modeling.}

\begin{figure}[h]
     \begin{center}
       \includegraphics[width=2cm]{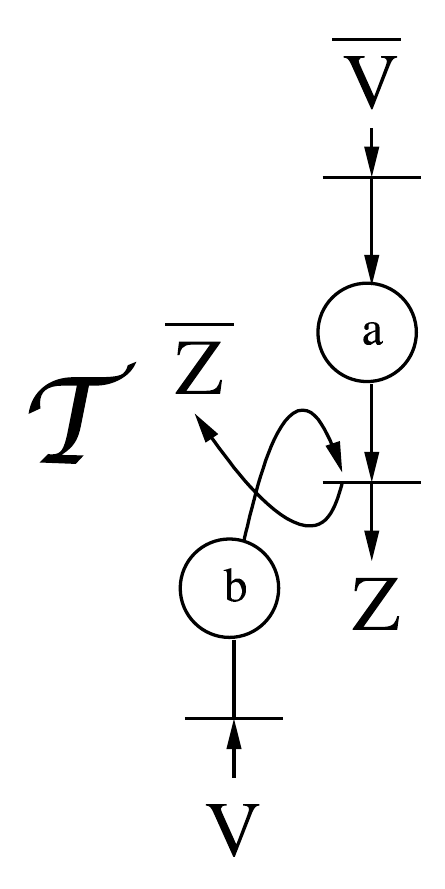}\hspace{5mm}
       \includegraphics[width=5cm]{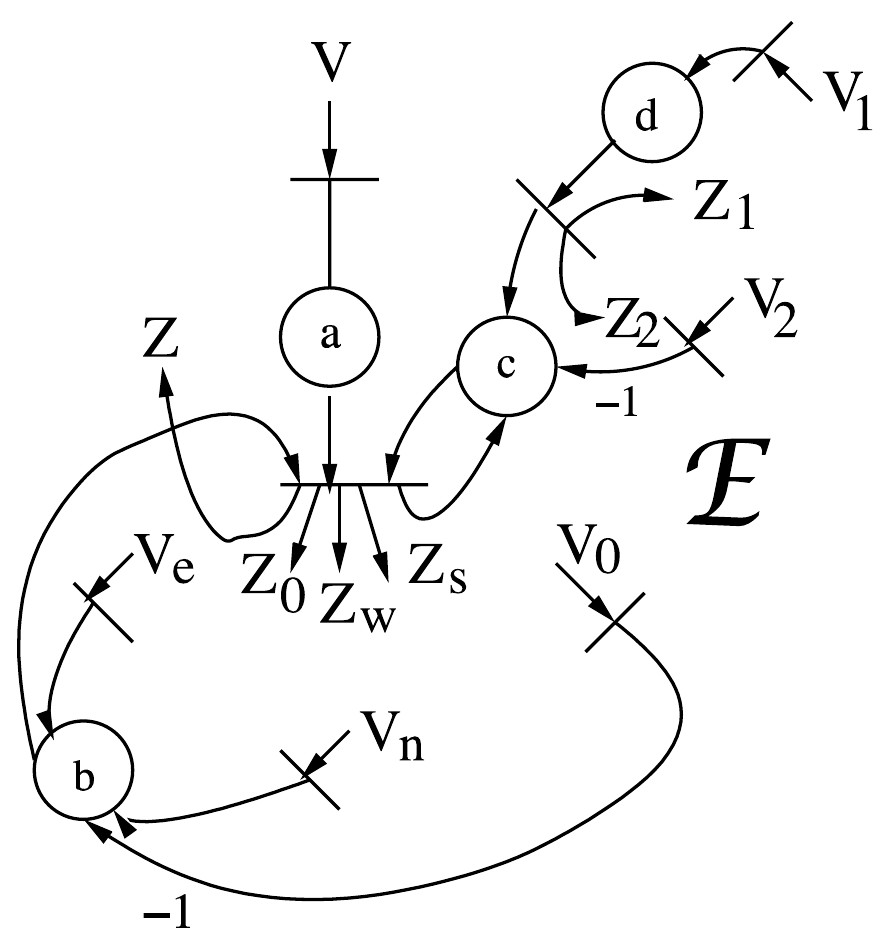}\hspace{1mm}
       \includegraphics[width=4.2cm]{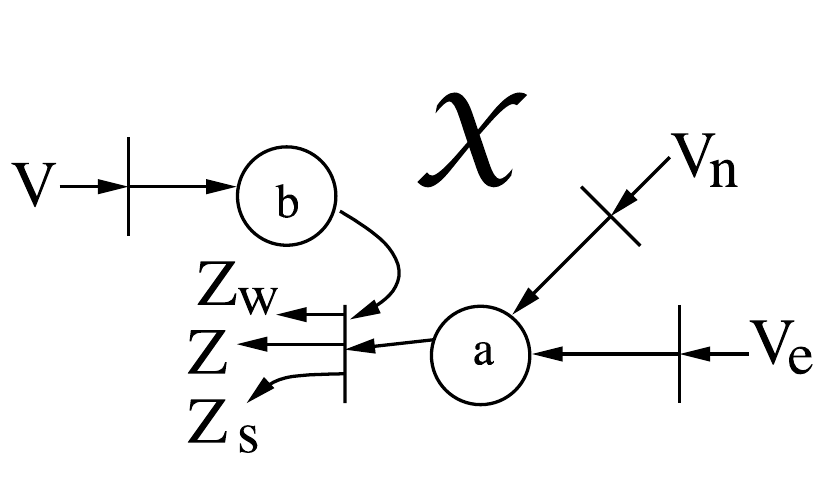}
       \caption{The three elementary Petri nets with which it is possible
       to build the Petri net of a regular set of roads on a torus by composition.}
       \label{elements}
     \end{center}
\end{figure}

\begin{figure}
  \begin{center}
    \includegraphics[width=14cm]{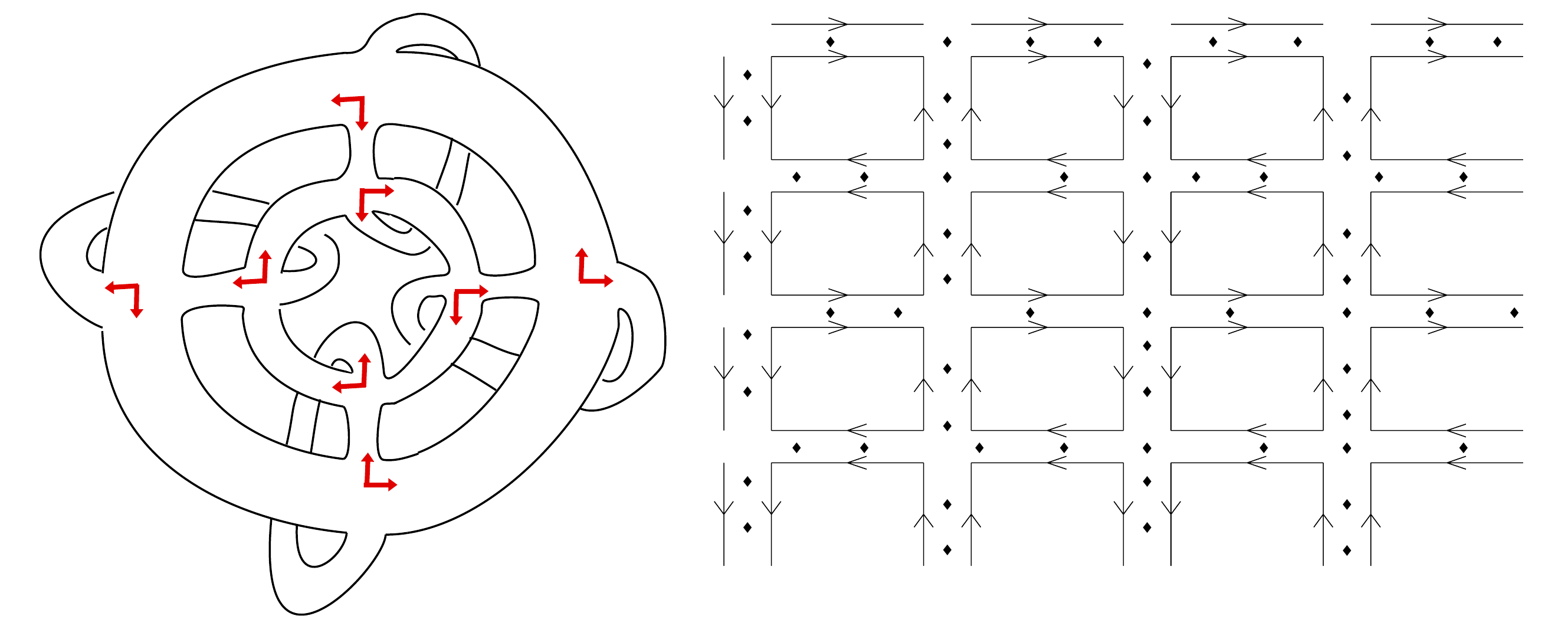}
    \caption{Roads on a torus of $4\times 2$ streets with its authorized turn at junctions (left) and tht asymptotic car distribution in the streets on a torus of $4\times 4$ streets obtained by simulation.}
    \label{ggreat}
  \end{center}
\end{figure}

We can generalize the modeling approach used in the case of one
junction to  derive the  fundamental diagram   of    a   regular
town on a torus    described  in Figure-\ref{ggreat} (left side). We
build the model by combining three elementary systems
(representing one section (T) , a junction input (E), and a junction
output (X)) with the composition operators described in the
generalized system theory subsection of the minplus section. The
details can be found in \cite{FARTh}.

The asymptotic vehicle distribution for a small town composed of
two North-South, South-North, East-West and West-East avenues is
given in Figure-\ref{ggreat} (right side). The fundamental diagram
presents a four-phase shape analogous to the case of two roads
with one junction. In this more general case, the role of the
non-priority road is played by a circuit of non-priority roads
which blocks the whole system when the circuit is full.

\section{Conclusion}
The Petri net modeling of traffic in a regular town can be done  thanks to
the introduction  of  negative  weights  on output  transition
edges.  The dynamics have a nice degree one homogeneous minplus
property but are no longer monotone. This loss of monotonicity
implies that the eigenvalue and the growth rate are no longer
equal. Experimental results show that they are close in the case
of two roads with one crossing where we are able to solve
explicitly the eigenvalue problem and compute numerically the
growth rate. The fundamental traffic diagram, which gives the
dependence of the average car-flow (given by the growth rate) on
the density, presents four phases that have traffic
interpretations.

This set of 1-homogeneous minplus systems seems to be a good class
of systems that we can  describe by two matrices, one in the
standard algebra and one in the minplus algebra. The  standard
compositions of these systems  can be easily described in terms of
these two matrices. These compositions of simple systems are useful
to build dynamics of large systems like that of regular town
traffic.

In a forthcoming and more traffic-oriented paper, we will describe
further the traffic interpretation of the four phases valid also
for more general systems like regular towns. The influence on the fundamental
diagram of the traffic control, using signal lights, will also be studied.

\section{Appendix}
\begin{lemma}\label{Reduced}
The eigenvalue problem~(\ref{Eigel}) of size $N$ can be reduced to the eigenvalue
problem~(\ref{RedEig}) of size 4. \end{lemma}
\begin{proof}
First let us verify that the solutions of (\ref{Eigel}) verify $\lambda \leq 1/4$.
Indeed (\ref{Eigel}) imply that $\lambda q_{n+m}\leq \bar{a}_{n+m}q_1q_{n+1}/(\lambda q_{n})$,
$\lambda q_{1}\leq a_{n}\sqrt{q_{n}q_{n+m}}$ and $\lambda q_{n+1}\leq a_{n+m}\sqrt{q_{n}q_{n+m}}$.
Multiplying these three inequalities we obtain $\lambda^4\leq \bar{a}_{n+m}a_{n}a_{n+m}=1$.

To obtain the result we have to eliminate $q_i,i\neq 1,n,n+1$,$n+m$, that is to solve
a linear minplus system which has a unique solution as soon as $\lambda < 1/2$. Indeed
the loops of the precedence graph associated to the linear system has all its loops positive
when $\lambda\leq 1/2$.

Moreover we can compute explicitly its solution. For  $i=2,n-1$~ we have~:
$$q_i=q_1\left[\bigotimes_{j=1}^{i-1}(a_j/\lambda)\right]\oplus \left[\bigotimes_{j=i}^{n-1}(\bar{a}_j/\lambda)\right]q_n\;,$$
for $i=n+1,n+m-1$~:
$$q_i=q_{n+1}\left[\bigotimes_{j=n+1}^{i-1}(a_j/\lambda)\right]\oplus \left[\bigotimes_{j=i}^{n+m-1}(\bar{a}_j/\lambda)\right]q_{n+m}\; .$$
Using this explicit solution in the four last equations of (\ref{Eigel}) we obtain the reduced system~(\ref{RedEig}).
\end{proof}

To  verify the results given Table~\ref{Array}, let us rewrite the
system~(\ref{RedEig}) with simplified notations~: $U=q_n$,
$V=q_{n+m}$, $X=q_1$, $Y=q_{n+1}$, $g=b_n$, $h=b_m$, $k=a_n$,
$l=a_{n+m}$,  w$n'=n-1$, and $m'=m-1$.
\begin{equation}  \begin{cases}
     U= \bar{k}XY / \lambda V \oplus gX/\lambda^{n'}\;,\\
     V=\bar{k}XY/\lambda^2 U\oplus hY/\lambda^{m'}\;,\\
     X=k\sqrt{UV}/\lambda \oplus \bar{g}U/\lambda^{n'}\;,\\
     Y=l\sqrt{UV}/\lambda \oplus \bar{h}V/\lambda^{m'} \;,
  \end{cases}\label{EigNew} \end{equation}
  with $\bar{k}=1-k-l\geq 0$, $\bar{g}=n'-g\geq 0$ and $\bar{h}=m'-h\geq 0$ are the free places in the
  crossing and the two roads.

\begin{lemma}
The eigenvector $(U,V,X,Y)$ of the minplus nonlinear system~(\ref{EigNew}) is given (using minplus notations) by~:
{\small \begin{equation}\begin{array}{|c||c|c|c|c|}
\hline
  & 0\leq d \leq \alpha & \alpha \leq d \leq \beta
&\min(\beta,\gamma) < d < \max(\beta,\gamma) & \gamma\leq d \leq 1 \\
 \hline \hline
 U& g/\lambda^{n'}  & g/\lambda^{n'}  & g/\lambda^{n'} & \bar{k}\bar{h} \\
 \hline
V & \lambda^{n'+2}/k^2g =hl/k\lambda^{m'}& 1\lambda^{n'-2}/k^2g=\lambda^{n'+2}/k^2g & \lambda^{n'+2}/k^2g & e/k^2\bar{k}\bar{h} \\
 \hline
 X& e  & e  & e & e \\
  \hline
Y& l/k & l/k & \lambda^{2+n'-m'}\bar{h}/k^2g= \lambda^4l/1k & e/k^2\bar{k}\\
\hline
\end{array}
\label{array2} \end{equation}}
where the eigenvalue $\lambda$ of the minplus nonlinear system~(\ref{EigNew}) has been given in
Theorem~\ref{maintheo} $\square$
\end{lemma}


\begin{proof}
Let us consider the four density regions corresponding to the four columns of this table~:
\begin{enumerate}
\item $0\leq d\leq\alpha$~: the first column of Table~\ref{Array} is solution of the standard linear
system~:
\begin{equation}   U= gX/\lambda^{n'},\quad
     V= hY/\lambda^{m'}, \quad
     X=k\sqrt{UV}/\lambda\;, \quad
     Y=l\sqrt{UV}/\lambda  \;,\label{sys1}\end{equation}
which is itself a solution of System~(\ref{EigNew}) since~:
\begin{enumerate}
\item  $\bar{k}XY / \lambda VU=\bar{k}kl/\lambda^3 =1/\lambda^3\geq 0$ (since $\lambda\leq 1/4$
indeed, once more, using~(\ref{EigNew}) we have $VXY\leq \bar{k}klVXY/\lambda^4$),
\item  $\bar{k}XY / \lambda^2 VU=1/\lambda^4 \geq 0$ ,
\item  $\bar{g}U/X\lambda^{n'}= (1/\lambda^2)^{n'}\geq 0$ ,
\item  $\bar{h}V/Y\lambda^{m'}=\bar{h}h/\lambda^{m'}=(1/\lambda^2)^{m'}\geq 0$ .
\end{enumerate}

Moreover, multiplying the 4 equalities of (\ref{sys1}) we obtain
$\lambda^{n'+m'+2}=ghlk$ which gives the value of $\lambda$ given in
Table~\ref{Array}.

\item $\alpha \leq d \leq \beta$~: the second column of Table~\ref{Array} is solution of the standard linear
system~:
\begin{equation}  U= gX/\lambda^{n'},\quad
      V= \bar{k}XY/\lambda^2 U\;, \quad
      X=k\sqrt{UV}/\lambda\;, \quad
      Y=l\sqrt{UV}/\lambda  \;,\label{sys2}\end{equation}
which is itself solution of System~(\ref{EigNew}) since~:
\begin{enumerate}
\item  $\bar{k}XY / \lambda VU= 1/\lambda^3\geq 0 $ ,
\item  $hY/V\lambda^{m'}=hglk/\lambda^{n'+m'+2}=d^{m'+n'+1}/(mn)^{1/4}\geq 0$ since $d\geq \alpha$ ,
\item  $\bar{g}U/X\lambda^{n'}=\bar{g}g/\lambda^{2 n'}=(1/\lambda^2)^{n'}\geq 0 $ ,
\item  $\bar{h}V/Y\lambda^{m'}=m'\lambda^{n'+2-m'}/hkgl=m'(2n'/m')^{1/4}/d^{m'+n'+1}\geq 0$ since $d\leq \beta$ .
\end{enumerate}

Moreover using the equality giving two expressions for the value of $V$ in (\ref{array2}) we obtain $\lambda=1/4$.

\item $\min(\beta,\gamma) \leq d \leq \max(\beta,\gamma)$~: the third column of Table~\ref{Array} is solution of the standard linear
system~:
\begin{equation}   U= gX/\lambda^{n'},\quad
     V= \bar{c}XY/\lambda^2 U\;, \quad
     X=k\sqrt{UV}/\lambda\;, \quad
     Y=\bar{h}V/\lambda^{m'}  \;,\label{sys3}\end{equation}
which is itself solution of System~(\ref{EigNew}) since~:
\begin{enumerate}
\item  $\bar{k}XY / \lambda VU= \lambda\geq 0 $ ,
\item  $hY/V\lambda^{m'}= h\bar{h}/\lambda^{ 2m'}=(1/\lambda^2)^{m'}\geq 0 $ since $\lambda\leq 1/4$ ,
\item  $\bar{g}U/X\lambda^{n'}= \bar{g}g/\lambda^{2 n'}=(1/\lambda^2)^{n'}\geq 0$ ,
\item  $l\sqrt{UV}/\lambda Y= 1/\lambda^4 \geq 0$ .
\end{enumerate}

Moreover using the equality giving two expressions for the value of $Y$ in (\ref{array2}) we obtain $\lambda^{2+n'-m'}m'=klgh=(m'+n'+1)d$ which gives the value of $\lambda$ given in Table~\ref{Array}.

\item $\gamma \leq d \leq 1$~: the fourth column of Table~\ref{Array} is solution of the standard linear
system~:
\begin{equation}    U= \bar{k}XY/  \lambda V\;,\quad
      V= \bar{k}XY/\lambda^2 U\;, \quad
      X=e\sqrt{UV}/\lambda\;, \quad
      Y=\bar{h}V/\lambda^{m'} \;,\label{sys4}\end{equation}
which is itself solution of System~(\ref{EigNew}) since~:
\begin{enumerate}
\item  $gX/\lambda^{n'}U=g/\bar{k}\bar{h}=ghkl/m=d^{n+m-1}/m\geq 0$ (since $d\geq \gamma$) ,
\item  $hY/\lambda^{m'}V=h\bar{h}=m'\geq 0$ ,
\item  $\bar{g}U/X\lambda^{n'}= \bar{g}\bar{k}\bar{h}\geq 0$ ,
\item  $l\sqrt{UV}/\lambda Y=lk\bar{k}=1$ .
\end{enumerate}

Moreover the compatibility of first two equalities of (\ref{sys4})
implies that $\lambda=0$.
\end{enumerate}
\end{proof}

\end{document}